\documentclass[12pt,twoside]{amsart}
\input{includeNice3}
\usepackage{mabliautoref}
\usepackage{amssymb}
\usepackage{amscd}
\usepackage{csquotes}
\usepackage{hyperref}
\usepackage[a4paper,margin=1in]{geometry}  
\usepackage{CJKutf8}
\usepackage{soul}

\usepackage{tikz-cd}
\usepackage{caption}

\usepackage{url}
\RequirePackage{booktabs, multirow}
\RequirePackage{pgf}

\newcommand{\cO}{\mathcal{O}}
\newcommand{\cM}{\mathcal{M}}

\theoremstyle{plain}
\newtheorem{thmA}{Theorem}

\title[A strong counterexample to the lc Beauville--Bogomolov decomposition]
{A strong counterexample to the log canonical Beauville--Bogomolov decomposition}

\author[F. Bernasconi, S. Filipazzi, Zs. Patakfalvi, and N. Tsakanikas]{Fabio Bernasconi, Stefano Filipazzi, Zsolt Patakfalvi, and Nikolaos Tsakanikas, with an appendix by Niklas M\"uller}



\thanks{FB is supported by the grant PZ00P2-21610 from the Swiss National Science Foundation. SF and NT are partly supported by the ERC starting grant $\#804334$. ZsP is supported by the grant $\#$200020B/192035 from the Swiss National Science Foundation, and by the ERC Starting grant $\#804334$. NM is supported by the grant $\#412744520$ (GRK2553) from the German research foundation.}
 
\address{Dipartimento di Matematica “Guido
Castelnuovo”, SAPIENZA Università di Roma, Piazzale Aldo Moro 5, I-00185, Roma, Italy}
\email{fabio.bernasconi@uniroma1.it}

\address{Department of Mathematics, Duke University,
120 Science Drive,
117 Physics Building,
Campus Box 90320,
Durham, NC 27708-0320,
USA}
\email{stefano.filipazzi@duke.edu}

\address{\'Ecole Polytechnique F\'ed\'erale de Lausanne, Chair of Algebraic Geometry, C3 635 
	(B\^atiment MA), Station 8, 1015 Lausanne, Switzerland 
}
\email{zsolt.patakfalvi@epfl.ch}

\address{\'Ecole Polytechnique F\'ed\'erale de Lausanne, Chair of Algebraic Geometry, C3 595 
	(B\^atiment MA), Station 8, 1015 Lausanne, Switzerland
}
\email{nikolaos.tsakanikas@epfl.ch}

\address{\'Ecole Polytechnique F\'ed\'erale de Lausanne, Chair of Algebraic Geometry, C3 615 
	(B\^atiment MA), Station 8, 1015 Lausanne, Switzerland.}
\email{{\tt niklas.mueller@epfl.ch}}

\DeclareMathOperator{\proj}{pr}
\DeclareMathOperator{\sbs}{\mathbf{B}} 

\newcommand{\Q}{\mathbb{Q}}

\newcommand{\C}{\mathbb{C}}

\begin{document}

    \begin{abstract}
        For every $d \geq 4$, we construct a $d$-dimensional, log canonical, $K$-trivial variety with the property that two general fibers of its Albanese morphism are not birational.
        We also show that this  provides a strong counterexample to the Beauville--Bogomolov decomposition in the log canonical setting.
        The construction can be adapted to construct a smooth quasi-projective variety of logarithmic Kodaira dimension 0 whose quasi-Albanese morphism has maximal variation.
        On the positive side, we show that the Albanese morphism for log canonical pairs with nef anti-canonical class is a locally stable family of pairs.
    \end{abstract}

    \maketitle

    \tableofcontents

\section{Introduction}

Throughout this paper we work over the field $\mathbb C$ of complex numbers.

Smooth projective varieties with torsion canonical divisor are one of the fundamental classes of varieties studied in birational geometry.
The Beauville--Bogomolov decomposition provides a structure theorem for these varieties.
It asserts that a smooth $K$-trivial variety can be decomposed, possibly after an étale cover, into a product of abelian varieties, strict Calabi--Yau varieties, and irreducible holomorphic symplectic varieties \cite{Bog74, Beau83}.
From this one can deduce the statement that is sometimes called the weak Beauville--Bogomolov decomposition: the Albanese morphism of a variety with trivial canonical divisor is isotrivial. This was first established by Calabi \cite{Cal57}.
The isotriviality of the Albanese morphism has been extended by Cao \cite{Cao19} to the case of smooth projective varieties with nef anti-canonical divisor.
This result has later been applied by Cao and H\"{o}ring \cite{CH19} to establish a Beauville--Bogomolov type decomposition of their universal cover.

From the perspective of the classification theory of varieties, a decomposition theorem for smooth projective varieties with trivial canonical class is not sufficient.
One would need variants allowing mild singularities (e.g., klt or log canonical).
In fact, such varieties are one of the main building blocks of the final outputs of the Minimal Modeal Program (MMP), singular Fano varieties and canonically polarized varieties being the other two. 

Following this motivation, in the series of articles \cite{GKP16, Dru18, GGK19, HP19}, an analog of the Beauville--Bogomolov decomposition for projective varieties with klt singularities and numerically trivial canonical class has been obtained. 
This result has been partially extended to the case of klt pairs with nef anti-canonical class in \cite{CCM, PZ19, MW21}, and it can be summarized by the following structure theorem.

\begin{thmA}[Decomposition theorem for klt pairs with nef anti-canonical divisor]
\cite{GKP16, Dru18, GGK19, HP19,CCM, PZ19, MW21}
    \label{thm:structure_thm_nef_anticanonical}
    Let $(X, \Delta)$ be a projective klt pair such that $-(K_X+\Delta)$ is nef. The following statements hold:
    \begin{enumerate}
        \item \label{item: isotrivial} the Albanese morphism $\alb_X \colon (X, \Delta) \to \Alb_X$ is an isotrivial morphism for the pair $(X,\Delta)$ and has connected fibers;
  
        \item \label{itm:klt}  if $K_X+\Delta \equiv 0$, then there exists a finite quasi-\'etale cover $\gamma \colon Y \to X$ such that 
        $$ (Y, \gamma^* \Delta)=(F, \Delta_F) \times A \times \prod Y_i \times \prod Z_i, $$
        where $F$ is a rationally connected variety and $\Delta_F = ( \gamma^* \Delta ) |_F$, $A$ is an abelian variety, the $Y_i$ are irreducible Calabi--Yau varieties and the $Z_i$ are irreducible symplectic varieties; and
        
        \item \label{itm:smooth} if $X$ is smooth and $\Delta=0$, then the universal cover $Y$  of $X$ admits a decomposition 
        $$ Y=F \times \bC^q \times \prod Y_i \times \prod Z_i, $$
        where $F$, $Z_i$ and $Y_i$ are as above and smooth.
    \end{enumerate} 
\end{thmA}

Note that the decomposition of the universal cover of a klt pair with nef anti-canonical class is still an open problem; see \cite[Conjecture 1.5]{MW21}.
Note also that in point \autoref{itm:smooth} the splitting does not happen on a finite cover, as instead in point \autoref{itm:klt}, but only on the universal one.
This can be traced back to the fact that, for any polarization, the polarized automorphism group is finite in the $K$-trivial klt case, but it is infinite in general in the $-K_X$ nef case.

\medskip

Varieties and pairs with log canonical singularities form the largest class of varieties for which the MMP is expected to hold. It is thus a central question to decide whether Theorem \autoref{thm:structure_thm_nef_anticanonical} holds for log canonical pairs with nef anti-canonical class.
However, one could then wonder which point of Theorem \autoref{thm:structure_thm_nef_anticanonical} should generalize to the log canonical case.
As the polarized automorphism groups can be positive dimensional in the log canonical case, even in the $K$-trivial case (e.g., take $(\bP^1, \{0\}+ \{\infty\})$ with $\sO_{\bP^1}(1)$ as polarization), one does not expect point \autoref{itm:klt} to extend to the log canonical case. Hence, even in the $K$-trivial log canonical case, the best one could hope for is that a decomposition as in point \autoref{itm:smooth} of Theorem \autoref{thm:structure_thm_nef_anticanonical} holds on the universal cover.

Thus, in this article we are primarily interested in the following question:

\begin{question}
\label{qtn:main}
Given a log canonical variety (or pair) $X$ with $K_X \sim_{\bQ} 0 $, does the quasi-étale universal cover $Y$ of $X$ (which is just the normalization over $X$ of the universal cover of $X_{\reg}$) admit a decomposition 
$$ Y=F \times \bC^q \times \prod Y_i \times \prod Z_i, $$
where the factors $F$, $Y_i$ and $Z_i$ are the singular versions of the factors  of point \autoref{itm:smooth} of Theorem \ref{thm:structure_thm_nef_anticanonical}?
\end{question}

Our main theorem is the following:

\begin{theorem}
    \label{thm:main}
    The answer to \autoref{qtn:main} is negative.
\end{theorem}

In fact, \autoref{thm:main} is an immediate corollary of the following more precise theorem, which states that the Albanese morphism fails, in general, to be isotrivial.

\begin{theorem} 
    \label{thm: counterexample_lc}
    For every integer $d \geq 4$, there exists a projective log canonical variety $X$ of dimension $d$ such that the following hold, where $\alb_X : X \to \Alb_X$ is the Albanese morphism of $X$:
    \begin{enumerate}
        \item $K_X \sim 0$;
        
        \item $\Alb_X$ is an elliptic curve;
        
        \item \label{item: not_bir_fiber} every fiber of $\alb_X$ is birational to exactly finitely many other fibers; and
    
        \item \label{item: univ_cover} the natural map $\pi_1(X_{\rm reg}) \to \pi_1(\Alb_X)$ is an isomorphism.
    \end{enumerate}
    Therefore, any quasi-\'etale cover of $X$ is induced by an \'etale cover of $\Alb_X$.
    In particular, the universal cover of $X_{\rm reg}$ admits a fibration to $\mathbb C$ where any fiber is birational to exactly countably many other fibers.
\end{theorem}

In the case of pairs, we also construct a counterexample to the Beauville--Bogomolov decomposition where the pair has plt singularities, see \autoref{prop: plt-counterexample}.
We note that the pair obtained in \autoref{prop: plt-counterexample} is crepant birational to the variety described in \autoref{thm: counterexample_lc}. The proofs of these two statements are given in \autoref{section:counterexamples_BB}.

We note that there exists an earlier example \cite[Example 6.3]{EIM23}, whose Albanese morphism is isotrivial, but does not split after a finite étale base change.
This shows that there is no log canonical Beauville--Bogomolov decomposition as in point \autoref{itm:klt} of Theorem \ref{thm:structure_thm_nef_anticanonical}.
However, this example really only exploits the fact that the polarized automorphism groups are positive dimensional.
In particular, the fibration does become split when one passes to the quasi-étale universal cover, and hence it does not give a counterexample to \autoref{qtn:main}.

We also treat the case of open varieties.
By \cite{Kaw81, Fuj24}, the quasi-Albanese morphism of a quasi-projective variety with logarithmic Kodaira dimension 0 is a dominant morphism with irreducible general fibers.
We show that this is optimal.

\begin{theorem} 
    \label{thm: non-isotriviality-quasialbanese}
    For every integer $d \geq 4$, there exists a smooth quasi-projective variety $U$ of dimension $d$ and logarithmic Kodaira dimension $\overline{\kappa}(U)=0$ such that the following hold:
    \begin{enumerate}
        \item the quasi-Albanese morphism $\alb_U \colon U \to \mathbb{G}_m$ is flat with irreducible fibers; and
        
        \item  every fiber of $\alb_U$ is birational to exactly finitely many other fibers.
    \end{enumerate}
\end{theorem}

We remark that the variety in \autoref{thm: counterexample_lc} is obtained by base change of a suitable compactification of the variety in \autoref{thm: non-isotriviality-quasialbanese}.%
\medskip

Despite the above counterexamples, we can show that the Albanese morphism of log Calabi--Yau pairs satisfies various desirable properties.
For example, it is a flat surjective morphism with semi-log canonical fibers.
This answers affirmatively Demailly--Peternell--Schneider's conjecture \cite[Conjecture 2]{DPS93} in the projective log canonical case.
Note that part \autoref{CZ_item} of our following result was already obtained in \cite{ChenZhang13}, even though it was explicitly stated only for smooth varieties.

\begin{theorem}
    \label{thm: lc_alb_contraction}
    Let $(X, \Delta)$ be a projective log canonical pair such that $-(K_X+\Delta)$ is nef. Then  $\alb_X \colon (X,\Delta) \to \Alb_X$ is a locally stable family of pairs. In particular, the following hold:
    \begin{enumerate}
        \item \label{CZ_item} $\alb_X \colon X \to \Alb_X$ is a contraction (i.e., $(\alb_X)_*\mathcal{O}_X=\mathcal{O}_{\Alb_X}$);
        
        \item \label{item: alb_red_fibers} $\alb_X$ is flat and has reduced fibers;
        
        \item every fiber $(F, \Delta|_F)$ of $\alb_X$ is a semi-log canonical pair; and

        \item \label{rmk:lcc_dominate} every log canonical center of $(X,\Delta)$ dominates $\Alb_X$.
    \end{enumerate}
\end{theorem}

Due to \autoref{thm: counterexample_lc}, we cannot expect $\alb_X \colon X \to \Alb_X$ to be isotrivial.
On the other hand, we show that the minimal log canonical centers of dlt log Calabi--Yau pairs are isotrivial onto their Albanese variety; see \autoref{thm:albanese_minimal_lcc} for the precise statement.

Lastly, one may wonder whether it is necessary for the pair in \autoref{thm: lc_alb_contraction} to be log canonical.
This is indeed the case: we construct a Fano $3$-fold with worse-than-log canonical singularities whose Albanese morphism is not surjective.

\begin{theorem}
    \label{thm: counterexample_non_lc}
    Let $C$ be a smooth projective curve of genus $g \geq 2$.
    Then there exists a normal projective 3-fold $X$ with $\Q$-Cartier and ample anti-canonical class ${-}K_X$ such that $\Alb_X \cong \Alb_C$ and the image of $\alb_X$ is identified with the image of $\alb_C$ under this isomorphism.
    In particular, $\alb_X$ is not surjective.
\end{theorem}

Similarly, \autoref{example:big_non-nef_nonsurjAlb} demonstrates that the assumption on the nefness of the anti-canonical class in \autoref{thm: lc_alb_contraction} cannot be dropped either.

\subsection*{Acknowledgments}

We would like to thank  H.\ Blum, S.\ Boissière, E.\ Floris, L.\ Ji, M.\ Mauri, S.-i.\ Matsumura, N.\ M\"uller, R.\ Svaldi, and C.\ Spicer for many fruitful conversations and comments on the topics of this work. The authors would also like to thank J.\ Moraga for helpful feedback on an earlier version of this work.

\section{Preliminaries}

\subsection{Notation}

\begin{enumerate}
    \item A variety $X$ is an integral separated scheme of finite type over $\mathbb C$.
    
    \item For a normal variety $X$, we denote by $\omega_X$ its canonical sheaf. 
    A canonical divisor $K_X$ is a Weil divisor on $X$ such that $\mathcal{O}_X(K_X)=\omega_X$.
    
    \item We say $(X, \Delta)$ is a \emph{pair} if $X$ is a normal variety and $\Delta$ is an effective $\mathbb{Q}$-divisor such that $K_X+\Delta$ is $\mathbb{Q}$-Cartier. For the definition of the singularities of the MMP we refer to \cite{kk-singbook}.

    \item We say that two pairs $(X_1,\Delta_1)$ and $(X_2,\Delta_2)$ are \emph{crepant birational} to one another if there exist proper birational morphisms $p_1 \colon Y \to X_1$ and $p_2 \colon Y \to X_2$ from a normal variety $Y$ such that
    $$
    p_1^*(K_{X_1}+\Delta_1)=p_2^*(K_{X_2}+\Delta_2) .
    $$

    \item Let $f \colon X \to Y$ be a proper morphism between varieties.
    We say that $f$ is a \emph{contraction} if $f_* \mathcal{O}_X=\mathcal{O}_Y$ holds.
    We say that $f$ is a \emph{fibration} if it is a contraction and moreover $\dim Y < \dim X$ holds.

    \item Let $f \colon X \to Y$ be a projective morphism of varieties.
    We say that $f$ is \emph{isotrivial} if it is locally trivial (i.e., isomorphic to the product family) in the \'etale topology.
    If we further assume that $X$ is endowed with the structure of a pair $(X,\Delta)$, then we say that $f \colon (X,\Delta) \to Y$ is \emph{isotrivial for the pair} $(X,\Delta)$ if it is locally trivial in the \'etale topology as a morphism of pairs.

    \item Let $f \colon X \to Y$ be a contraction of quasi-projective varieties.
    We say that $f$ is \emph{birationally isotrivial} if, for any two closed points $y_1, y_2 \in Y$, the varieties $X_{y_1}$ and $X_{y_2}$ are birational to each other.
    Furthermore, if $X$ is endowed with the structure of a pair $(X,\Delta)$, then we say that $f \colon (X,\Delta) \to Y$ is \emph{crepant birationally isotrivial} if, for any two closed points $y_1, y_2 \in Y$, the pairs $(X_{y_1},\Delta|_{X_{y_1}})$ and $(X_{y_2},\Delta|_{X_{y_2}})$ are crepant birational to one another.
    Lastly, we say that $f$ (resp.\ $f \colon (X,\Delta) \to Y$) is \emph{generically birationally isotrivial} (resp.\ \emph{generically crepant birationally isotrivial}) if it is birationally isotrivial (resp.\ crepant birationally isotrivial) over a non-empty open subset $V \subset Y$.

    \item Given a smooth quasi-projective variety $U$, we define the \emph{logarithmic Kodaira dimension} $\overline{\kappa}(U)$ of $U$ to be the Kodaira dimension $\kappa(X,K_X+D)$, where $X$ is a smooth projective variety and $D$ is an snc divisor on $X$ such that $U = X \setminus D$.
    
    \item A variety $X$ is said to be \emph{uniruled} if there exists a variety $Y$ of dimension $\dim Y = \dim X - 1$ and a dominant rational map $Y \times \bP^1 \dashrightarrow X $. We say that $X$ is \emph{rationally connected} (resp.\ \emph{rationally chain connected}) if any two general points on $X$ can be joined by a rational curve (resp.\ a chain of rational curves).
\end{enumerate}

\subsection{Generalized pairs}

We recall here the notion of a generalized pair, which was originally introduced by Birkar and Zhang \cite{BZ16}, as well as the usual classes of singularities of generalized pairs. For details about the language of b-divisors we refer to \cite[\S~2.1]{FS23}.

\begin{definition} 
    \label{def: gen_pair}
    A \emph{generalized sub-pair} $(X, \Delta, \bM )/Z$ over a variety $Z$ is the datum of 
    \begin{enumerate}
        \item a normal quasi-projective variety $X$, endowed with a projective morphism $X \to Z$;
        \item a $\mathbb{Q}$-divisor $\Delta$ on $X$; and
        \item a $\mathbb{Q}$-b-Cartier $\mathbb{Q}$-b-divisor $\bM $ that is b-nef over $Z$,
    \end{enumerate}
    satisfying that $K_X+\Delta+\bM_X$ is a $\mathbb{Q}$-Cartier $\mathbb{Q}$-divisor.
    If $\Delta$ is effective, then we say that $(X,\Delta,\bM) / Z$ is a \emph{generalized pair} over $Z$.
    If $Z=\Spec(\mathbb C)$, we omit $Z$ from the notation and simply write $(X,\Delta,\bM)$.
    Moreover, if $\bM=0$, then we also drop $\bM$ and the word \enquote{generalized} from the notation, and we retrieve thus the usual notions of \emph{sub-pair} and \emph{pair}.
\end{definition}

\begin{definition}
    Let $(X, \Delta, \bM)/Z$ be a generalized pair and let $\pi \colon X' \to X$ be a projective birational morphism from a normal variety $X'$. We define a $\bQ$-divisor $\Delta'$ on $X'$ via the identity
    \[
    K_{X'} +\Delta' +\bM_{X'} =\pi^*(K_X +\Delta+\bM_{X}).
    \]
    Given a prime divisor $E$ on $X'$, we define the \emph{generalized log discrepancy} of $E$ with respect to $(X,B,\bM)/Z$ to be 
    \[ 
    a_E(X, \Delta, \bM) \coloneqq 1-\text{coeff}_E(\Delta') . 
    \]
\end{definition}

\begin{definition}
    \label{dfn: gen_singularities}
    We say that a generalized pair $(X, \Delta, \bM )/Z$ is \emph{generalized log canonical} (resp.\ \emph{generalized klt}) if $a_E(X, \Delta, \bM ) \geq 0$ (resp.\ $a_E(X, \Delta, \bM ) >0$) for every divisor $E$ over $X$.
\end{definition}

\subsection{Generalized log Calabi--Yau pairs}

\begin{definition}
    \label{dfn: gen_log_CY}
    A projective generalized pair $(X,B,\mathbf M)$ is called \emph{generalized log Calabi--Yau} if it is generalized log canonical and $K_X +B +\mathbf M _X \sim_{\mathbb Q} 0$ holds.
    If $\mathbf M =0$, then we say that $(X,B)$ is a \emph{log Calabi--Yau pair}.
\end{definition}

\begin{example}
    Let $(X,\Delta)$ be a projective log canonical pair such that $-(K_X+\Delta)$ is nef.
    Define the b-divisor $\mathbf M \coloneqq \overline{-(K_X+\Delta)}$.
    Then $(X,\Delta,\mathbf M)$ is a generalized log Calabi--Yau pair by \cite[Remark 4.2(4)]{BZ16}.
\end{example}

\begin{lemma} 
    \label{lemma:canonical}
    Let $(X,\Delta,\mathbf M)$ be a generalized log Calabi--Yau pair.
    Then $X$ is not uniruled if and only if $X$ has canonical singularities and $K_X \sim_\Q 0$, $\Delta=0$, and $\mathbf M _X \sim_\Q 0$ hold.
\end{lemma}

\begin{proof}
    Assume first that $X$ is not uniruled.
    Let $h \colon (X',\Delta',\mathbf M) \to (X,\Delta,\mathbf M)$ be a $\Q$-factorial generalized dlt modification of $(X,\Delta,\mathbf M)$ and note that $ K_{X'} + \Delta' + M_{X'} \sim_\Q 0 $; see \cite[Theorem 2.9]{FS23}.
    Since $X'$ is not uniruled either, it follows from \cite[Lemma 2.18]{TX23b}, which relies crucially on \cite[Corollary 0.3]{BDPP}, that $X'$ has canonical singularities, and also that $K_{X'} \sim_\Q 0$, $\Delta'=0$, and $\mathbf M _{X'} \equiv 0$ hold; in particular, we have $\mathbf M _{X'} \sim_\Q 0$. Thus, by construction, we deduce that $\Delta=0$, $\mathbf M _X = h_* \mathbf M_{X'} \sim_\Q 0$, and that $X$ has canonical singularities.
    In particular, we have $K_X \sim_\Q 0$.
    Finally, the converse follows readily from \cite[Corollary 0.3]{BDPP}.
\end{proof}

The following statement is a generalization of the results of Ambro \cite{Amb04,Amb05} to the setting of generalized pairs. It was obtained in \cite{Fil20,FS23}, and we include it here for the reader's convenience.

\begin{theorem}[Adjunction and inversion of adjunction for fiber spaces] 
    \label{thm: adj_fiber_space}
    Let $(X, \Delta, \bM)/Z$ be a generalized pair and let $f \colon X \to Y$ be a contraction of normal varieties over $Z$.
    Suppose $K_X + \Delta +\bM _X \sim_{\mathbb{Q},f} 0$ and that the generic fiber is generalized log canonical.
    Then there exists a generalized pair structure $(Y, B_Y, \mathbf{N})/Z$ on $Y$ such that 
    \[ 
    K_X + \Delta +\bM _X \sim_{\mathbb{Q}} f^*(K_Y+ B_Y +\mathbf{N}_Y) .
    \]
    Moreover, adjunction and inversion of adjunction hold:
    $(X, \Delta, \bM)/Z$ is generalized log canonical if and only if so is $(Y, B_Y, \mathbf{N})/Z$.
    Furthermore, if $(X, \Delta, \bM)/Z$ is generalized klt, then so is $(Y, B_Y, \mathbf{N})/Z$.
\end{theorem}

\begin{proof}
    In the case $Z=\Spec(\mathbb C)$, the first assertion is \cite[Theorem 1.4]{Fil20}, while for a general base $Z$ it is proven in \cite[Theorem 2.20]{FS23}.
    The adjunction statement is proven in \cite[Proposition 4.16]{Fil20}; note that op.\ cit.\ is stated for $Z=\Spec(\mathbb C)$, but the proof is independent of what $Z$ is.
\end{proof}

\begin{lemma} 
    \label{lem:surj_map_av}
    Let $(X,\Delta,\mathbf M)$ be a generalized log Calabi--Yau pair. Let $f \colon X \to A$ be a surjective morphism to an abelian variety and consider its Stein factorization
    $ X \longrightarrow Z \overset{g}{\longrightarrow} A$.
    Then $Z$ is an abelian variety and $g$ is an \'etale cover.
\end{lemma}

\begin{proof}
    By \autoref{thm: adj_fiber_space}, $Z$ inherits from $(X,\Delta,\mathbf M)$ the structure of a generalized log Calabi--Yau pair; we denote it by $(Z,\Delta_Z,\mathbf N)$.
    Since $Z$ admits a finite surjective morphism to an abelian variety, it is not uniruled, and it follows now from \autoref{lemma:canonical} that $Z$ has canonical singularities and that $K_Z \sim_{\mathbb Q} 0$.

    Since $A$ is an abelian variety, we may choose a nowhere vanishing section $\sigma$ of $\omega_A$.
    Its pull-back $g^*\sigma$ is a non-zero section of $\omega_Z$.
    Since $K_Z\sim_{\mathbb Q}0$ holds, it follows that $g^*\sigma$ is nowhere vanishing as well.
    In particular, we have $K_Z \sim 0$ and $g$ has to be unramified in codimension 1, since otherwise $g^*\sigma$ would vanish along some divisor. By purity of the branch locus \cite[Tag 0BMB]{stacks-project}, we conclude that $g$ is \'etale. Therefore, $Z$ itself is an abelian variety by \cite[\S~18, Theorem (Serre--Lang)]{MumfordAV}.
\end{proof}

\subsection{Albanese morphism}

The \emph{Albanese morphism} $\alb_X \colon X \to \Alb_X$ of a normal projective variety $X$ has been constructed by Serre \cite{Ser60}. It holds that $\Alb_X \simeq (\Pic^0_X)_{\text{red}}^\vee$ by \cite[Remark 9.5.25]{FGA}. In particular, the \emph{Albanese variety} $\Alb_X$ of $X$ is an abelian variety. 

We collect below some properties of the Albanese morphism that we will need in this paper.

\begin{lemma}
    \label{lem:RCC_surj_Alb}
    Let $X$ be a proper normal variety. 
    If $X$ is rationally chain connected, then the Albanese variety $\Alb_X$ of $X$ is trivial and, in particular, the Albanese morphism $\alb_X \colon X \to \Alb_X$ of $X$ is surjective.
\end{lemma}

\begin{proof}
    As $\Alb_X$ is an abelian variety and every morphism from a rational curve to an abelian variety is necessarily constant, we conclude.
\end{proof}

\begin{lemma}
    \label{lem:Alb_RCC_fibers}
    Let $f \colon X \to Y$ be a contraction between proper normal varieties.
    If every fiber of $f$ is rationally chain connected, then the induced map $\Alb(f) \colon \Alb_X \to \Alb_Y$ is an isomorphism.
\end{lemma}

\begin{proof}
    By the universal property of the Albanese morphism, we obtain the following commutative diagram:
    \begin{center}
        \begin{tikzcd}[row sep = large, column sep = large]
            X \arrow[d, "\alb_X" swap] \arrow[r, "f"] & Y \arrow[d, "\alb_Y"] \\
            \Alb_X \arrow[r, "\Alb(f)" swap] & \Alb_Y.
        \end{tikzcd}
    \end{center}
    If $F$ is a closed fiber of $f$, then it is rationally chain connected by assumption. Thus, it is contracted to a closed point of $\Alb_X$, since abelian varieties contain no rational curves. By the Rigidity Lemma \cite[Lemma 1.15(b)]{Deb01}, we infer that $\alb_X$ factors through $f$,
    that is, there exists a morphism $\varphi \colon Y \to \Alb_X$ such that $\varphi \circ f = \alb_X$: 
    \begin{center}
        \begin{tikzcd}[row sep = large, column sep = large]
            X \arrow[d, "\alb_X" swap] \arrow[r, "f"] & Y \arrow[d, "\alb_Y"] \arrow[dl, "\exists \varphi" swap] \\
            \Alb_X \arrow[r, "\Alb(f)" swap] & \Alb_Y .
        \end{tikzcd}
    \end{center}
    By using the universal property of the Albanese morphism, we first verify that $\alb_Y = \Alb(f) \circ \varphi$, and then we easily conclude that $\Alb_X \simeq \Alb_Y$.
\end{proof}

\begin{remark}
    If the morphism $f$ is equidimensional, then it is sufficient to ask that the generic fiber of $f$ is rationally chain connected, as rational chain connectedness specializes in equidimensional families by \cite[Corollary IV.3.5.2]{Kol96}.
\end{remark}

\section{The Albanese morphism for log canonical pairs with nef anti-canonical divisor}

In this section we prove several structure results for the Albanese morphism of pairs $(X,\Delta)$ whose anti-canonical class $-(K_X+\Delta)$ has some positivity property. We begin with the proof of \autoref{thm: lc_alb_contraction}.

\begin{proof}[Proof of \autoref{thm: lc_alb_contraction}]
    We first show that $\alb_X \colon X \to \Alb_X$ is a contraction. To this end, set $\mathbf{M} \coloneqq \overline{{-}(K_X+\Delta)}$ and consider the generalized log Calabi--Yau pair $(X,\Delta, \mathbf{M})$.
    Let $V$ be the image of $\alb_X \colon X \to \Alb_X$ and
    let $X \to W \to V$ be the Stein factorization.
    By \autoref{thm: adj_fiber_space} there exist $\Delta_W$ and $\mathbf{N}$ such that $(W, \Delta_W, \mathbf{N})$ is a generalized log Calabi--Yau pair.
    As $W$ admits a generically finite map to an abelian variety, it is not uniruled, so \autoref{lemma:canonical} implies that $K_W \sim_\Q 0$, $\Delta_W=0$, $\mathbf{N}_W \sim_\Q 0$, and that $W$ has canonical singularities.
    In particular, $\mathbf N$ descends on $W$.

    Consider now the normalization $V^n \to V$ and the induced map $W \to V^{n} \to V$.
    By construction, $W \to V^n$ is finite, so by \cite[Lemma 1.1]{FG12} we may find a boundary $B$ on $V^n$ such that $(V^n, B)$ is a klt log Calabi--Yau pair.
    Since $V^n$ is not uniruled, \autoref{lemma:canonical} implies that $K_{V^n} \sim_\Q 0$, $B=0$ and $V^n$ has canonical singularities.
    Since $V$ generates $\Alb_X$, we can apply \cite[Corollary 3.5]{Uen73} to deduce that $V=\Alb_X$, and thus $W \to V$ is an \'etale covering of abelian varieties by \autoref{lem:surj_map_av}.
    By the universal property of the Albanese morphism, we conclude that $W=V=\Alb_X$, which shows item \autoref{CZ_item} of the theorem.

    Finally, let $D$ be a reduced divisor in $\Alb_X$ with simple normal crossing support. Then the pair $(\Alb_X,D)$ is log canonical, so the generalized pair $(\Alb_X,D,\mathbf N)$ is in turn generalized log canonical, since $\mathbf N$ descends on $\Alb_X = W$ by the property $\mathbf{N}_{\Alb_X} \sim_\Q 0$.
    By \autoref{thm: adj_fiber_space}, $ \big( X,\Delta+ (\alb_X)^*D,\mathbf M \big)$ is also generalized log canonical,
    and hence $\big( X,\Delta+ (\alb_X)^*D \big)$ is log canonical, since $\mathbf M$ descends on $X$ by definition.
    By \cite[Corollary 4.55]{Kol23}, $\alb_X \colon (X,\Delta) \to \Alb_X$ is a locally stable family.
    In turn, items (b) and (c) are formal consequences of this fact, see \cite[Definition-Theorem 4.7]{Kol23}, whereas item \autoref{rmk:lcc_dominate} follows from \cite[Corollary 4.56]{Kol23}.
\end{proof}

The next example shows that the fibers of the Albanese morphism for log Calabi--Yau pairs are not necessarily irreducible.
Therefore, one cannot expect isotriviality without any further assumption.

\begin{example} 
    \label{ex: 2dim}
    Consider the log Calabi--Yau surface pair 
    \[ (X, \Delta) \coloneqq \big( \mathbb{P}^1 \times E, \, \left\{ 0 \right\} \times E + \left\{ \infty \right\} \times E \big) , \]
    where $E$ is an elliptic curve. Denote by $\pi_E$ the projection $\mathbb P^1 \times E \to E$.
    Let $f \colon Y \to X$ be the blow-up of $X$ at a closed point 
    $x \in \left\{ 0 \right\} \times E$ and let $(Y, \Delta_Y)$ be the log Calabi--Yau pair crepant birational to $(X,\Delta)$.
    Then the Albanese morphism $Y \to E$ is the composite map $\pi_E \circ f$, which is not isotrivial, since it has a non-irreducible fiber.
\end{example}

We now study the Albanese morphism of the log canonical centers of log canonical pairs with nef anti-canonical class.

\begin{lemma} \label{lem:alb_lcc_surjective}
    Let $(X, \Delta)$ be a projective log canonical pair such that $-(K_X+\Delta)$ is nef.
   Let $V$ be a log canonical center of $(X,\Delta)$ and let $V^n$ be its normalization.
   Then the natural morphism $\alb(\iota) \colon \Alb_{V^n} \to \Alb_X$ induced by $\iota \colon V^n \to X$ is surjective.
\end{lemma}

\begin{proof}
    By items \autoref{CZ_item} and \autoref{rmk:lcc_dominate} in \autoref{thm: lc_alb_contraction}, the induced morphism $\alb_X \circ \iota \colon V^n \to \Alb_X$ is surjective.
    By the universal property of the Albanese morphism, this morphism factors through $\Alb_{V^n}$, 
    and it follows that $\alb(\iota) \colon \Alb_{V^n} \to \Alb_X$ is surjective, as claimed.
\end{proof}

\begin{theorem}
    \label{thm:albanese_minimal_lcc}
    Let $(X, \Delta)$ be a projective dlt pair such that $-(K_X+\Delta)$ is nef.
    Let $V$ be a minimal log canonical center of $(X,\Delta)$ and let $(V,\Delta_V)$ be the klt pair induced by adjunction of $(X,\Delta)$ to $V$.
    Then the following hold:
    \begin{enumerate}
        \item\label{part a thm albanese minimal lcc} $\alb_V \colon (V,\Delta_V) \to \Alb_V$ is isotrivial for the pair $(V,\Delta_V)$;
        
        \item\label{part b thm albanese minimal lcc} the natural morphism $\Alb_V \to \Alb_X$ is surjective; and
        
        \item\label{part c thm albanese minimal lcc} the morphism $(V,\Delta_V) \to A$ is isotrivial for the pair $(V,\Delta_V)$, where $A$ is the abelian variety obtained as the Stein factorization of $\Alb_V \to \Alb_X$.
    \end{enumerate}
\end{theorem}

\begin{proof}
    Since $(X,\Delta)$ is dlt, $V$ is normal by \cite[Theorem 4.16]{kk-singbook}.
    The pair $(V, \Delta_V)$ is induced by $(X, \Delta)$ via higher Poincar\'e residue maps \cite[4.18.4]{kk-singbook} and it satisfies $(K_X+\Delta)|_V \sim_{\mathbb{Q}} K_V+\Delta_V$.
    Since $V$ is a minimal log canonical center, $(V,\Delta_V)$ is klt by \cite[Theorem 4.19]{kk-singbook}.
    Moreover, since $K_V+\Delta_V$ is obtained by adjunction of $K_X+\Delta$, $-(K_V+\Delta_V)$ is nef.
    Thus, part \autoref{part a thm albanese minimal lcc} of the statement follows from part \autoref{item: isotrivial} of Theorem \ref{thm:structure_thm_nef_anticanonical}.

    Part \autoref{part b thm albanese minimal lcc} follows immediately from \autoref{lem:alb_lcc_surjective}.

    Since $\Alb_V \to \Alb_X$ is surjective, $A$ is the Stein factorization of a surjective morphism between abelian varieties.
    By \autoref{lem:surj_map_av}, $A$ is an abelian variety itself.
    Then, part \autoref{part c thm albanese minimal lcc} follows from \cite[Theorem 4.1]{MW21}.
\end{proof}

\begin{remark}
    We observe that, in part \autoref{part c thm albanese minimal lcc} of \autoref{thm:albanese_minimal_lcc}, we actually obtain a stronger result, namely that $(V,\Delta_V) \to A$ is a locally constant fibration with respect to the pair $(V,\Delta_V)$ in the sense of \cite[Definition 2.3]{MW21}.
\end{remark}

The following two examples demonstrate that item \autoref{part b thm albanese minimal lcc} of \autoref{thm:albanese_minimal_lcc} is optimal.

\begin{example}
    \label{example:1_Thm3.3_optimal}
    Consider the log Calabi--Yau pair $(\mathbb{P}^2, E)$, where $E$ is an elliptic curve.
    Then $E$ is the minimal log canonical center of $(\mathbb{P}^2, E)$ and $\Alb_E=E$, while $\Alb_{\mathbb{P}^2}=\Spec(\mathbb C)$.
\end{example}

\begin{example}
    \label{example:2_Thm3.3_optimal}
    Let $E$ be an elliptic curve and consider the log smooth log Calabi--Yau surface pair 
    \[ \big( \mathbb{P}^1 \times E, \, \left\{ 0 \right\} \times E + \left\{ \infty \right\} \times E \big) . \]
    Let $P$ be a 2-torsion point on $E$ and consider the following $(\mathbb Z / 2 \mathbb Z)$-action on $\mathbb P ^1 \times E$:
    \begin{equation*}
        \begin{split}
            \tau \colon \mathbb P ^1 \times E &\to \mathbb P ^1 \times E\\
            (x,y) &\mapsto (x^{-1},y+P).
        \end{split}
    \end{equation*}
    Observe that $\tau$ has no fixed points and preserves the boundary, and hence it induces an automorphism of the pair 
    \[ \big( \mathbb P ^1 \times E, \{0\} \times E + \{\infty\} \times E \big) . \]
    Let $(X,\Delta)$ be the pair obtained as the quotient.
    Since $\tau$ has no fixed points, the quotient map $\mathbb P ^1 \times E \to X$ is \'etale and $X$ is smooth.
    Since $\tau$ interchanges the two components of $\{0\} \times E + \{\infty\} \times E$, the boundary $\Delta$ is irreducible and isomorphic to $E$.
    In particular, $(X,\Delta)$ is again log smooth with numerically trivial log canonical class.
    Furthermore, since $\tau$ is Galois, 
    $H^0 \big(X, \mathcal{O}_X(m(K_X+\Delta)) \big)$ coincides with the $(\mathbb{Z}/2\mathbb{Z})$-invariant part of
    $$
    H^0 \Big(\mathbb P ^1 \times E, \mathcal{O}_{\mathbb P ^1 \times E} \big( m(K_{\mathbb{P}^1\times E}+\{0\} \times E + \{\infty\} \times E) \big) \Big),
    $$
    which is generated by $(\frac{dx}{x} \wedge dy)^{\otimes m}$.
    This shows that $K_X+\Delta$ is a non-trivial 2-torsion divisor.
    Moreover, we have an induced morphism $X \to E'$, where $E'$ is the elliptic curve obtained by quotienting $E$ by the action $y \mapsto y+P$.
    The morphism $X \to E'$ is a $\mathbb P ^1$-fibration, so it coincides with $\alb_X$ by \autoref{lem:Alb_RCC_fibers}.
    In turn, the natural map $\Alb_\Delta \to \Alb_X$, which coincides with the morphism $\Delta \to E'$, is an \'etale 2:1 cover, so its fibers are not connected.
\end{example}

We conclude this section by showing that one can relate the Albanese morphism of a pair with big anti-canonical class with the Albanese morphism of a canonically defined subvariety of $X$. Before giving the precise statement, we need to recall the following two definitions.

Fix a normal projective variety $X$. The stable base locus of a $\Q$-Cartier $\Q$-divisor $\Gamma$ on $X$ is denoted by $\sbs(\Gamma)$. Given a $\Q$-Cartier $\Q$-divisor $D$ on $X$, the \emph{non-nef locus} of $D$ (also referred to as the \emph{diminished base locus} of $D$ in the literature) is defined as 
\[ \nnef(X,D) \coloneqq \bigcup_{m \geq 1} \sbs \big( D+ \frac{1}{m}A \big) , \] 
where $A$ is an ample divisor on $X$. Note that the above definition does not depend on the ample divisor chosen and that $\nnef(X,D) = \emptyset$ if and only if $D$ is nef. Now, given a subset $V$ of $X$, we say that $X$ is \emph{rationally connected modulo $V$} if either (a) $V = \emptyset$ and $X$ is rationally connected, or (b) $V \neq \emptyset$ and there exists an irreducible component $Z$ of $V$ such that for any general point $x \in X$ there exists a rational curve passing through $x$ and intersecting $Z$.

According to \cite[Theorem 1.2]{BP11}, if $(X, \Delta)$ is a projective pair such that $-(K_X+\Delta)$ is big, then there exists an irreducible component $V$ of $\nnef(-K_X-\Delta) \cup \nklt(X, \Delta)$ such that $X$ is rationally connected modulo $V$. With this result in mind, we are now ready to state and prove the promised statement. Specifically:

\begin{corollary} 
    \label{lem:induced_map_alb}
    Let $(X, \Delta)$ be a projective pair such that $-(K_X+\Delta)$ is big. Consider the locus $V$ defined above 
    and denote by $V^n$ its normalization.
    Then the induced map $\Alb_{V^n} \to \Alb_X$ is surjective.
\end{corollary}

\begin{proof}
    Consider the commutative diagram given by the universal property of the Albanese morphism:
    \begin{center}
        \begin{tikzcd}[row sep = large, column sep = large]
            V^n \arrow[d, "\alb_{V^n}" swap] \arrow[r, "\sigma"] & X \arrow[d, "\alb_X"] \\
            \Alb_{V^n} \arrow[r, "\Alb(\sigma)" swap] & \Alb_X.
        \end{tikzcd}
    \end{center}
    Let $x \in \Alb_X$. 
    Since abelian varieties do not contain rational curves, 
    it follows from \cite[Theorem 1.2]{BP11} that the fiber $\alb_X^{-1}(x)$ intersects $V$.
    In particular, $\alb_X^{-1}(x)$ intersects the image of $\sigma$.
    This implies that the image of $\Alb(\sigma) \colon \Alb_{V^n} \to \Alb_X$ contains $\alb_X(X) \subseteq \Alb_X$.
    Since $\Alb(\sigma)$ is a group homomorphism and $\alb_X(X)$ generates $\Alb_X$, we deduce the desired surjectivity of $\Alb(\sigma)$.
\end{proof}

\section{Counterexamples to the log canonical Beauville--Bogomolov decomposition}
\label{section:counterexamples_BB}

In this section we construct several counterexamples to the Beauville--Bogomolov decomposition in the log canonical case. We also prove \autoref{thm: counterexample_lc} and \autoref{thm: non-isotriviality-quasialbanese}.

\subsection{A crepant birationally isotrivial example}

As a warm-up, we construct an example of a 3-dimensional log Calabi--Yau pair whose Albanese morphism is not isotrivial.

\begin{proposition} 
    \label{ex: 3-dimensional}
    There exists a log smooth 3-dimensional log Calabi--Yau pair $(\mathcal{Y}, \mathcal{B})$ such that 
    \begin{enumerate}
        \item $\alb_\mathcal{Y} \colon \mathcal{Y} \to E$ is a smooth fibration onto an elliptic curve;
        
        \item all fibers are weak del Pezzo surfaces of degree 4; and
        
        \item two general fibers $\mathcal{Y}_t$ and $\mathcal{Y}_s$ are not isomorphic.
    \end{enumerate}
\end{proposition}

\begin{proof}
    Consider the product $\mathbb P ^2 \times E$, where $E$ is an elliptic curve embedded as a smooth cubic inside $\mathbb{P}^2$. We may then consider $E \times E \subset \mathbb P^2 \times E$ as a divisor.
    We observe that $(\mathbb P ^2 \times E, E \times E)$ is a log Calabi--Yau pair.
    Furthermore, by construction, the projection $\mathbb P ^2 \times E \to E$ is the Albanese morphism of $\mathbb P ^2 \times E$.
    
    We fix four distinct closed points $P_1,P_2,P_3,P_4 \in E \subset \mathbb P ^2$. We denote by $\pi \colon \mathcal{X} \to \mathbb P ^2 \times E$ the blow-up of $\mathbb P ^2 \times E$ along the disjoint union of curves $\bigcup_{i=1}^4 \left\{P_i \right\} \times E$, and by $\mathcal{E}$ the strict transform of $E \times E$ on $\mathcal{X}$.
    By the inclusion $\{P_1,P_2,P_3,P_4\} \times E \subset E \times E$, 
    we obtain
    $$ 
        K_{\mathcal{X}}+\mathcal{E} = \pi^* \big( K_{\mathbb P ^2 \times E}+(E \times E) \big) ,
    $$ 
    and we also infer that $\mathcal{E} \to E \times E$ is an isomorphism.
    Let $\Delta$ be the preimage of the diagonal in $E \times E$ under said isomorphism, and let $\psi \colon \mathcal{Y} \to \mathcal{X}$ denote the blow-up of $\mathcal{X}$ along $\Delta$.
    As above, we have 
    $$ 
        K_{\mathcal{Y}}+\mathcal{B} = \psi^*(K_{\mathcal{X}}+\mathcal{E}) ,
    $$ 
    where $\mathcal{B}$ denotes the strict transform of $\mathcal{E}$.
    Furthermore, $\mathcal{B} \to \mathcal{E}$ is again an isomorphism.
    
    To summarize, we have obtained a log Calabi--Yau pair $(\mathcal{Y},\mathcal{B})$ with Albanese morphism $f \colon \mathcal{Y} \to E$.
    By construction, every fiber of $f$ is a blow-up of $\mathbb P ^2$ at 5 points on the given elliptic curve $E$.
    For a general fiber, these points are all distinct, while for 4 fibers, 2 of these points are infinitely close.
    In particular, all fibers of $f$ are smooth weak del Pezzo surfaces, and, with the exception of 4 closed fibers, they are actually del Pezzo surfaces of degree 4.
    
    We are left to verify that, for two general points $P$ and $Q$ in $E$, the closed fibers $X_P$ and $X_Q$ are not isomorphic.
    Indeed, let $\psi\colon X_P \to X_Q$ be an isomorphism.
    As $X_P$ is the blow-up at $\mathbb{P}^2$ at the closed points $P_1, \ldots, P_4$ and $P$ with exceptional divisors $E_i$ over $P_i$, the images $F_i \coloneqq \psi(E_i)$ and $F \coloneqq \psi(E)$ are disjoint $(-1)$-curves on $X_Q$ and we can blow them down by $\pi \colon X_Q \to \mathbb{P}^2$.
    Therefore, the automorphism $\psi$ is a lifting of the automorphism $\varphi\colon \mathbb{P}^2 \to \mathbb{P}^2$ sending $P_1, \dots, P_4$ and $P$ to $Q_1, \dots Q_5$. 
    Up to the action of $\mathrm{PGL}_3(k)$ on the codomain, we can assume that $P_i=Q_i$ for $i=1, \dots, 4$. 
    This implies that $\varphi(P)=Q_5$, which contradicts the fact that $P$ and $Q$ are general.
\end{proof}

The counterexamples to the Beauville--Bogomolov decomposition presented in \autoref{ex: 2dim} and \autoref{ex: 3-dimensional} are nevertheless isotrivial after a birational modification, so they are not sufficient to prove \autoref{thm: counterexample_lc}.

\subsection{A plt counterexample to the Beauville--Bogomolov decomposition} 
\label{sec:key_example}

In this subsection we construct, for every integer $d \geq 4$, a plt log Calabi--Yau $d$-fold whose Albanese morphism fails to be birationally isotrivial; see \autoref{prop: plt-counterexample}.

In this example we are interested in pencils of degree $d$ hypersurfaces inside $\mathbb{P}^d$.
We choose a family of degree $d \geq 4$ to construct the counterexample because, thanks to \cite[Theorem 30]{Kol19} and \cite[Main Theorem 2]{Kol19b}, we have direct tools to show that two general elements are not birational.

Now, let $X=V(f) \subset \mathbb{P}^d$ be a general smooth hypersurface of degree $d$.
In particular, $V(f,x_d)$ is a smooth $K$-trivial variety $S$.
We want to perturb $f$ while preserving the intersection $V(f,x_d)$.

The space of homogeneous polynomials of degree $(d-1)$ in $x_0,\ldots,x_d$ has dimension $\binom{2d-1}{d}$,
while the space of homogeneous polynomials of degree $d$ in $x_0,\ldots,x_d$ has dimension $\binom{2d}{d}$.
We consider the affine morphism
\[
\begin{split}
\iota \colon \mathbb A ^{\binom{2d-1}{d}} &\rightarrow \mathbb A ^{\binom{2d}{d}}\\
q &\mapsto f+x_dq
\end{split}
\]
where $q$ denotes a homogeneous polynomial of degree $d-1$.
The morphism is injective and therefore corresponds to a subvariety of dimension $\binom{2d-1}{d}$.
Since $f$ is general, it is not divisible by $x_d$.
In particular, the image of $\iota$ does not contain the origin of $\mathbb A ^{\binom{2d}{d}}$.
Then the rational map $\mathbb A^{\binom{2d}{d}} \dashrightarrow \mathbb P^{\binom{2d}{d}-1}$ is defined along $\im(\iota)$.
Since $\im(\iota)$ is closed and does not contain the origin, it does not contain any orbit under the natural $\mathbb G _m$-action on $\mathbb A^{\binom{2d}{d}}$.
It follows that the image of $\mathbb A ^{\binom{2d-1}{d}}$ in $\mathbb P ^{\binom{2d}{d}}$ is an irreducible locally closed subset of dimension $\binom{2d-1}{d}$.

Since $V(f)$ is smooth and being smooth is an open condition, we deduce that for a general choice of $q \in \mathbb A ^{\binom{2d-1}{d}}$, $V(f+x_dq)$ is smooth.
Note that $\dim \mathrm{Aut}(\mathbb P ^d)=(d+1)^2-1$ and the image of $\mathbb A ^{\binom{2d-1}{d}}$ in $\mathbb P ^{\binom{2d}{d}-1}$ has dimension $\binom{2d-1}{d}$.
Since $d \geq 4$, we have
$$ 
    \binom{2d-1}{d}>(d+1)^2-1 .
$$
Thus, for a general choice of two homogeneous polynomials $q_1$ and $q_2$ of degree $(d-1)$, the hypersurfaces $V(f+x_dq_1)$ and $V(f+x_dq_2)$ are smooth and not projectively equivalent, and hence not birational.
Since $q_1$ and $q_2$ are general, we may also assume that $V(q_1-q_2)$ and $V(q_1-q_2,x_d)$ are smooth.

\begin{notation} \label{notation:ex}
To sum up,
\begin{enumerate}
    \item[(i)] $V(f)$ is a general smooth Fano hypersurface in $\mathbb{P}^d$ of degree $d$ such that $V(f,x_d)$ is a smooth $K$-trivial hypersurface of degree $d$ in $\mathbb{P}^{d-1}$;
    
    \item[(ii)] $q_1$ and $q_2$ are general homogeneous polynomials of degree $(d-1)$ such that the varieties $V(f+x_dq_1), V(f+x_dq_2), V(q_1-q_2)$ and $V(q_1-q_2,x_d)$ are smooth;
    
    \item[(iii)] we consider the pencil of degree $d$ hypersurfaces 
    \begin{equation} 
        \label{eq: pencil_quadrics}
        \mathcal{X} \coloneqq V \big( y_0(f+x_dq_1)+y_1(f+x_dq_2) \big) \subset \mathbb P ^d_{[x_0:\ldots :x_d]} \times \mathbb P ^1_{[y_0:y_1]},
    \end{equation}
    together with the projection $\pi \colon \mathcal{X} \to \mathbb{P}^1$;
    
    \item[(iv)] the fibers $\mathcal{X}_{[0:1]}$ and $\mathcal{X}_{[1:0]}$ of $\pi$ are smooth and not birational to one another.
    \end{enumerate}
\end{notation}

In the following claims we analyze the singularities and the geometry of $\mathcal{X}$.

\begin{claim} 
    \label{cl-birational_fibers}
    Any two general closed fibers of $\pi \colon \mathcal{X} \to \mathbb{P}^1$ are not birational.
\end{claim}

There is a non-empty open subset $U \subset \mathbb P ^1$ such that every fiber over $U$ is a smooth degree $d$ hypersurface of dimension $d-1$.
Notice that $[0:1],[1:0] \in U$. 
By \cite[Corollary 1.5]{Fuj19} every fiber over $U$ is K-stable.
By \cite[Theorem 1.1]{BX19} the moduli stack of K-stable $\Q$-Fano varieties is separated.
Since the fibers over $[0:1]$ and $[1:0]$ are not isomorphic and the moduli stack is separated, it follows that two general fibers over $U$ are not isomorphic.
In turn, by \cite[Theorem 30]{Kol19} and \cite[Main Theorem 2]{Kol19b}, two general fibers over $U$ are not birational. 

\begin{claim}
    Set $\mathcal{D} \coloneqq \mathcal{X} \cap V(x_d)$. Then the pair $(\mathcal{X}, \mathcal{D})$ is log canonical.
\end{claim}

Notice that $y_0(f+x_dq_1)+y_1(f+x_dq_2)$ has bidegree $(d,1)$ and that $\mathcal{X}$ is a Fano variety of dimension $d$.
By construction, since $\mathcal{D}$ is the restriction to $\mathcal{X}$ of a divisor of bidegree $(1,0)$, $\pi \colon (\mathcal{X},\mathcal{D}) \to \mathbb P ^1$ is an lc-trivial fibration; see \cite[Definition 2.16]{FS23}.
We have 
\begin{equation*}
  K_\mathcal{X}+\mathcal{D} = \big( K_{\mathbb P ^d \times \mathbb P ^1}+\mathcal{X}+V(x_d) \big) \big|_\mathcal{X} \sim \mathcal{O}_{\mathcal{X}}(0,-1) \sim \pi^*\mathcal{O}_{\mathbb{P}^1}(-1) ,
\end{equation*}
so the canonical bundle formula applied to $(\mathcal{X}, \mathcal{D}) \to \bP^1$ has the form
\[
K_{\mathcal{X}} + \mathcal{D} \sim \pi^*(K_{\mathbb P ^1}+B+\bM_{\mathbb{P}^1}),
\]
where $\deg(B+\bM_{\mathbb{P}^1})=1$.
Since $\bM_{\mathbb{P}^1}$ is pseudoeffective and $B$ is effective, we have $\deg(B) \leq 1$.
In particular, since the generic fiber of $(\mathcal{X},\mathcal{D}) \to \mathbb P ^1$ is log canonical (log smooth even), it follows from \autoref{thm: adj_fiber_space} that $(\mathcal{X},\mathcal{D})$ is log canonical.

By construction, $\mathcal{D}$ is given by the equation $\left\{(y_0+y_1)f=0 \right\} \subset \mathbb{P}^{d-1} \times \mathbb{P}^1$, where $\mathbb P ^{d-1}$ is identified with the hyperplane $V(x_d)$ in $\mathbb P ^d$.
Thus, $\mathcal{D}$ is the union of $\mathbb{P}^{d-1} \times \{[1:-1]\}$ and $S \times \mathbb{P}^1$.
By the definition of the boundary divisor in the canonical bundle formula (see \cite[\S~2]{Amb04}), $[1:-1]$ has coefficient at least 1 in $B$.
Since $B$ has degree at most 1, we conclude that $B=[1:-1]$.
Then, by difference, we obtain $\mathbf{M}_{\mathbb{P}^1}=0$.

\begin{claim} 
    \label{claim: canonical-sing}
    $\mathcal{X}$ has canonical singularities.
\end{claim}

By adjunction, $-K_{\mathcal{X}}$ is ample and Cartier.
Thus, to show that $\mathcal{X}$ is canonical, it suffices to show that it is klt.
To this end, since $-K_{\mathcal{X}}$ is ample, we may define a generalized log Calabi--Yau pair $(\mathcal{X},0,\mathbf M)$, where $\mathbf M \coloneqq \overline{-K_{\mathcal{X}}}$.
By \autoref{thm: adj_fiber_space}, we have an induced generalized pair $(\mathbb P ^1,\Delta,\mathbf N)$ satisfying $K_{\mathcal{X}}+\mathbf{M}_{\mathcal{X}} \sim \pi^*(K_{\mathbb{P}^1}+\Delta+\mathbf{N}).$
Since $\mathbf M$ descends on $\mathcal{X}$, we have $\Delta \leq B$.
In particular, by \autoref{thm: adj_fiber_space}, $\mathcal{X}$ is klt over $\mathbb P ^1 \setminus \{[1:-1]\}$.
To conclude, by \autoref{thm: adj_fiber_space} once more, it suffices to show that $\Delta < B =[1:-1]$.
Assume by contradiction that this is not the case.
Then there is a log canonical center of $\mathcal{X}$ contained in $\mathcal{X}_{[1:-1]}$.
By construction, we have $\mathcal{X}_{[1:-1]}=V(x_d(q_1-q_2)) \subset \mathbb P^d$, so
\[ 
\mathcal{X}_{[1:-1]} \cong \mathbb P ^d \cup V(q_1-q_2) .
\]
Since $V(q_1-q_2)$ is a smooth variety of dimension $d-1$ by assumption,
$\mathcal{X}_{[1:-1]}$ is smooth away from $V(x_d) \cap V(q_1-q_2) \subset \mathbb P ^d$.
Therefore, the log canonical center of $\mathcal{X}$ is contained in $V(x_d) \cap V(q_1-q_2)$.
Yet, $(\mathcal{X},\mathcal{D})$ is log canonical, and $\mathcal{D}$ contains $V(x_d) \cap V(q_1-q_2)$, reaching thus a contradiction.
In turn, we have $\Delta < B= [1:-1]$, and the claim follows.

\begin{claim}
    $\mathcal{X}$ is not $\Q$-factorial.
\end{claim}

To prove the claim, it suffices to exhibit two prime Weil divisors whose intersection has codimension 3 and not 2.
Denote by $T$ the irreducible component of $\mathcal{X}_{[1:-1]}$ given by the degree $(d-1)$ hypersurface $V(q_1-q_2) \subset \mathbb P^d$, and by $\mathcal{D}^h$ the horizontal part of $\mathcal{D}$, which is isomorphic to $S \times \mathbb P ^1$.
We have $\mathcal{D}^h \subset \mathbb P ^{d-1} \times \mathbb P ^1$, where $\mathbb P ^{d-1}$ is identified with $V(x_d) \subset \mathbb P ^{d}$.
Hence, $\mathcal{D}^h \cap \mathcal{X}_{[1:-1]}=S \subset \mathbb P^{d-1}$.
In turn, $\mathcal{D}^h \cap T$ corresponds to the intersection in $\mathbb P ^{d-1}=V(x_d) \subset \mathbb P ^d$ between the degree $d$ hypersurface $S$ and the degree $(d-1)$ hypersurface $V(q_1-q_2,x_d)$, which is a codimension 3 variety.

\medskip

The following result builds on the previous claims and constitutes the key step towards the construction of the desired counterexample to the log canonical Beauville--Bogomolov decomposition.

\begin{proposition} 
    \label{prop: first_approx}
    There exists a projective plt pair $(\tilde{\mathcal{X}}, \tilde{\mathcal{D}}^{h})$ of dimension $d$ such that the following hold:
    \begin{enumerate}
        \item there exists a contraction $\tilde{
        \pi} \colon \tilde{\mathcal{X}} \to \mathbb{P}^1$;
        
        \item $K_{\tilde{\mathcal{X}}} + \tilde{\mathcal{D}}^{h} \sim \tilde{\pi}^*K_{\mathbb{P}^1}$; and
        
        \item every fiber of $\tilde{\pi}$ is birational to exactly finitely many other fibers.
    \end{enumerate}
\end{proposition}

\begin{proof}
We follow the notation from the previous claims and their proofs.
Now, since $\mathcal{X}$ is canonical by \autoref{claim: canonical-sing}, we may consider a small $\mathbb Q$-factorialization $\hat{\mathcal{X}} \to \mathcal{X}$; see \cite[Corollary 1.37]{kk-singbook}.
We denote by $\hat{\mathcal{D}}$ the strict transform of $\mathcal{D}$, and similarly, we write $\hat T$ for the strict transform of $T$ and $\hat \pi$ for the contraction onto $\mathbb{P}^1$.
Next, we run an MMP to contract $\hat T$.

\begin{claim}
    For $0 < \epsilon \ll 1$, the pair $(\hat{ \mathcal{X}},\hat{\mathcal{D}}-\epsilon \hat{\mathcal{D}}^v)$ is plt.
\end{claim}

Indeed, by construction, $(\hat{\mathcal X},\hat{\mathcal D})$ is isomorphic to $(\mathcal{X},\mathcal{D})$ over a non-empty open subset of $\mathbb P ^1$.
In particular, $(\hat{\mathcal X},\hat{\mathcal D})$ is log smooth and plt over a non-empty open subset of $\mathbb P ^1$.
Thus, to check that $(\hat{ \mathcal{X}},\hat{\mathcal{D}}-\epsilon \hat{\mathcal{D}}^v)$ is plt, we need to show that it has no log canonical centers that are vertical over $\mathbb P ^1$.
Since $(\hat{ \mathcal{X}},\hat{\mathcal{D}})$ is log canonical, no log canonical center can be contained in $\hat{\mathcal{D}}^v$.
Then the generic point of any vertical log canonical center is contained in $\hat{T} \setminus \hat{\mathcal{D}}^v$.
Since $T$ is smooth and it is Cartier away from $T \cap \mathcal{D}^v$, it follows that $\mathcal{X} \setminus \mathcal{D}^v$ is smooth in a neighborhood of $T \setminus \mathcal{D}^v$.
Therefore, $\hat{\mathcal X} \setminus \hat{\mathcal{D}}^v \to \mathcal{X} \setminus \mathcal{D}^v$ is the identity in a neighborhood of $\hat T \setminus \hat{\mathcal{D}}^v$.
Thus, the generic point of any vertical log canonical center is contained in the smooth locus of $\hat{\mathcal{X}}$ and inside $\hat{T} \setminus \hat{\mathcal{D}}^v$.
Then, by the inclusion $\mathcal{D}^h \cap T \subset \mathcal{D}^v \cap T$, it follows that $\hat{\mathcal{D}}^h \cap (\hat{T} \setminus \hat{\mathcal{D}}^v)=\emptyset$.
But then, the generic point of any vertical log canonical center is disjoint from the boundary of $(\hat{ \mathcal{X}},\hat{\mathcal{D}}-\epsilon \hat{\mathcal{D}}^v)$ and inside the smooth locus of $\hat{ \mathcal{X}}$.
In turn, no such log canonical center can exist, and hence $(\hat{ \mathcal{X}},\hat{\mathcal{D}}-\epsilon \hat{\mathcal{D}}^v)$ is plt, as claimed.

\begin{claim}
    We can run a $(K_{\hat{ \mathcal{X}}}+\hat{\mathcal{D}}-\epsilon \hat{\mathcal{D}}^v)$-MMP over $\mathbb{P}^1$, $\hat{ \mathcal{X}} \dashrightarrow \tilde{\mathcal{X}}$, which contracts the divisor $\hat{T}$.
\end{claim}

Indeed, we have
\[
    K_{\hat{\mathcal{X}}}+\hat{\mathcal{D}}-\epsilon \hat{\mathcal{D}}^v\sim_{\mathbb Q,\hat \pi} -\epsilon \hat{\mathcal{D}}^v \sim_{\mathbb Q,\hat \pi} \epsilon \hat T,
\]
where $\hat \pi$ denotes the induced morphism from $\hat{\mathcal{X}}$ to $\mathbb{P}^1$.
Then, by \cite[Lemma 2.10]{Lai11}, $\hat T$ is in the relative diminished base locus of $(\hat{\mathcal{X}},\hat{\mathcal{D}}-\epsilon \hat{\mathcal{D}}^v)$.
Thus, after finitely many steps of a relative $(K_{\hat{\mathcal{X}}}+\hat{\mathcal{D}}-\epsilon \hat{\mathcal{D}}^v)$-MMP with scaling of an ample divisor over $\mathbb{P}^1$ (see \cite{HX13}), $\hat T$ is contracted.
We denote the resulting model by $\tilde{\mathcal{X}}$ and the strict transform of $\hat{\mathcal{D}}$ on $\tilde{\mathcal{X}}$ by $\tilde{\mathcal{D}}$.
Since $K_{\hat{\mathcal{X}}}+\hat{\mathcal{D}}$ is relatively trivial over $\mathbb P^1$ by construction, $(\mathcal{X},\mathcal{D})$ is crepant birational to $(\tilde{\mathcal{X}},\tilde{\mathcal{D}})$.
Then, by construction, $\tilde{\mathcal{X}}_{[1:-1]}$ is irreducible and we have $\tilde{\mathcal{D}}^v=\tilde{\mathcal{X}}_{[1:-1]}$.
\end{proof}

Using the example constructed in \autoref{prop: first_approx}, it is now immediate to construct a counterexample to the Beauville--Bogomolov decomposition in the log setting.

\begin{proposition} 
    \label{prop: plt-counterexample}
    For every integer $d \geq 4$ there exists a projective plt pair $(\mathcal{Y}, \mathcal{B})$ of dimension $d$ such that the following hold:
    \begin{enumerate}
        \item $K_{\mathcal{Y}} + \mathcal{B} \sim 0$;
        
        \item $\Alb_{\mathcal{Y}}$ is an elliptic curve; and
        
        \item every fiber of $\alb_{\mathcal{Y}}$ is birational to exactly finitely many other fibers.
    \end{enumerate}
\end{proposition}

\begin{proof}
    If we consider the plt pair $(\tilde{\mathcal{X}},\tilde{\mathcal{D}}^h)$ constructed in \autoref{prop: first_approx}, then we have
    \[
        K_{\tilde{\mathcal{X}}/\mathbb{P}^1}+\tilde{\mathcal{D}}^h \sim 0.
    \]
    In particular, by \autoref{thm: adj_fiber_space}, $(\tilde{\mathcal{X}},\tilde{\mathcal{D}}^h) \to \mathbb P^1$ is a relatively minimal locally stable family of pairs; see \cite[\S~2.1]{Kol23}.
    Since the base change of a locally stable family remains a locally stable family, by base change to an elliptic curve $E \to \mathbb{P}^1$, we obtain a relatively minimal locally stable family $({\mathcal{Y}},\mathcal{B}) \to E$.
    By the compatibility of the relative log canonical divisor of a locally stable family under base change, we obtain
    \[
        0 \sim K_{\mathcal{Y}/E}+\mathcal{B} \sim K_{\mathcal{Y}}+\mathcal{B} .
    \]
    By construction, observe that the pair $(\mathcal{Y}, \mathcal{B})$ is plt and that the morphism $({\mathcal{Y}},\mathcal{B}) \to E$, which corresponds to the Albanese fibration of $\mathcal{Y}$, is not crepant birationally isotrivial, since any two general fibers are not even birational.
\end{proof}

\begin{remark}
    We observe that in \autoref{prop: plt-counterexample}, consistently with \autoref{thm:albanese_minimal_lcc}, the family induced by adjunction to the minimal log canonical center $\mathcal{B} \to \Alb_{\mathcal{Y}}$ is isotrivial.
\end{remark}

\begin{remark}
    We here explain why the proof of the klt case in \cite[Appendix A]{PZ19} breaks down in the above example. 
    A crucial ingredient in the proof of \cite[Proposition A.11]{PZ19} is that, in the klt case, if a relatively ample divisor is pseudoeffective, then it is actually nef. 
    In our case, $\mathcal{B}$ is effective and relatively ample, but not nef, as we explain below.
    
    First, na\"ively, if $\mathcal{B}$ were nef, then we could then apply the klt case \cite{PZ19} to the pair $\big( \mathcal{Y},(1-\epsilon)\mathcal{B} \big)$, for some $0< \epsilon < 1$.
    We now verify that $\mathcal{B}$ is not nef in geometric terms.
    As $\mathcal{B}$ is the pull-back of $\tilde{\mathcal{D}}^h$ by a finite surjective map, it is sufficient to show that $\tilde{\mathcal{D}}^h$ is not nef. To this end,
    fix $P \in S \setminus V(q_1-q_2,x_d)$ and consider the rational curve  $\Sigma \coloneqq \{P\} \times \mathbb P ^1 \subset \mathbb P ^d \times \mathbb P ^1$ contained in $S \times \mathbb{P}^1$.
    The divisor $\mathcal{D}^{h}$ is Cartier in a neighborhood of $\Sigma$, and a direct computation shows that $\mathcal{D}^{h} \cdot \Sigma =-1$.
    Since $\tilde{\mathcal{X}}$ is isomorphic to $\mathcal{X}$ in a neighborhood of $\Sigma$, it follows that $\tilde{\mathcal{D}}^{h}$ is not nef either.
\end{remark}

\subsection{The universal cover does not split}

We show that, if we endow the variety $\mathcal{Y}$ constructed in \autoref{prop: plt-counterexample} with the structure of a complex analytic space, then the Albanese morphism $\alb_{\mathcal{Y}}$ does not split after passing to the universal cover of the Albanese variety $\Alb_{\mathcal{Y}}$.

\begin{proposition} 
    \label{prop: univ_cover_boundary}
    Let $(\mathcal{Y}, \mathcal{B})$ be the pair constructed in \autoref{prop: plt-counterexample}.
    Consider the fiber product:
    \begin{center}
        \begin{tikzcd}[row sep = large, column sep = large]
            (\mathcal{Z}, \mathcal{B}_{\mathcal Z}) \arrow[r, "\pi"] \arrow[d, "h" swap] & (\mathcal{Y}, \mathcal{B}) \arrow[d, "\alb_{\mathcal{Y}}"] \\
            \mathbb{C} \arrow[r] & \Alb_{\mathcal{Y}} ,
        \end{tikzcd}
    \end{center}
    where $\mathbb{C} \to \Alb_{\mathcal{Y}}$ is the universal cover.
    Then
    \begin{enumerate}
        \item $\pi_1(\mathcal{Z}_{\reg} \setminus \mathcal{B}_{\mathcal{Z}})=1$; and
        
        \item two general fibers of $h \colon (\mathcal{Z}, \mathcal{B}_{\mathcal Z}) \to \mathbb{C}$ are not birational to each other.
    \end{enumerate}
\end{proposition}

Thus, even the copy of $\mathbb C$ that uniformizes $E \cong \Alb_\mathcal{Y}$ cannot be split off after taking the universal cover of $\mathcal{Y}_{\rm reg} \setminus \mathcal{B}$.

\begin{proof}
    Denote by $(Y,B)$ a general fiber of $({\mathcal{Y}},\mathcal{B}) \to E$.
    By construction, $Y$ is a smooth degree $d$ hypersurface of dimension $(d-1)$ and $B$ is a smooth $K$-trivial variety obtained as a hyperplane section.
    Since $d \geq 4$, by the Lefschetz hyperplane theorem for homotopy groups \cite[Theorem 3.1.21]{Laz} we have $\pi_1 (Y)=1$ and 
    \[ 
        H_2(Y,\mathbb Z)=H_2(\mathbb P ^d,\mathbb Z)=\mathbb Z .
    \]
    By \cite[Corollary 2.2]{Lib06}, $\pi_1(Y \setminus B)$ coincides with the cokernel of the map
    \begin{equation}
        \label{kronecker}
        \begin{split}
            H_2(Y,\mathbb Z) &\to \mathbb Z \\
            \alpha &\mapsto \langle \alpha , [B] \rangle ,
        \end{split}
    \end{equation}
    where $[B]$ denotes the class of $B$ in $H^2(X,\mathbb Z)$ and $\langle -,-\rangle$ denotes the Kronecker pairing.
    Recall now that $B$ is the restriction to $Y$ of a hyperplane in $\mathbb P ^d$.
    Then, by the identification $H_2(Y,\mathbb Z)=H_2(\mathbb P ^d,\mathbb Z)$, we may find a class $\alpha \in H_2(Y,\mathbb Z)$ such that $\langle \alpha , [B] \rangle = 1$.
    Thus, the map in \autoref{kronecker} is an isomorphism, and hence $\pi_1(Y \setminus B)=1$. 
    Since all fibers of $\alb_{\mathcal{Y}}$ are reduced (either by construction or by item \autoref{item: alb_red_fibers} of \autoref{thm: lc_alb_contraction}), by \cite[Lemma 1.5]{Nor83} we have the exact sequence of fundamental groups
    \[
        \pi_1(Y \setminus B) \to \pi_1(\mathcal{Y}_{\rm reg}\setminus \mathcal{B}) \to \pi_1(E) \to 1.
    \]
    Thus, $\pi_1(\mathcal{Y}_{\rm reg}\setminus \mathcal{B}) \to \pi_1(E)$ is an isomorphism.
    From this, we infer that $\pi_1(\mathcal{Z}_{\reg} \setminus \mathcal{B}_{\mathcal{Z}})=1$, and by \autoref{prop: plt-counterexample} two general fibers of $h$ are not birational to each other.
\end{proof}

Therefore, a na\"ive Beauville--Bogomolov decomposition as in \cite{CH19, CCM} cannot be expected for log canonical log Calabi--Yau pairs with integral coefficients (notice that this example is plt).

\subsection{Quasi-Albanese morphism}

In this section we show that the quasi-Albanese morphism of an open variety of logarithmic Kodaira dimension 0 is not generically birationally isotrivial. For the theory of quasi-Albanese varieties, morphisms, and their universal properties we refer to the recent article of Fujino \cite{Fuj24}.


\begin{proof}[Proof of \autoref{thm: non-isotriviality-quasialbanese}]
We utilize the notation in \autoref{sec:key_example}.
By construction, $(\tilde{\mathcal{X}},\tilde{\mathcal{D}}^h)$ is a plt pair with the property that
\[
K_{\tilde{\mathcal{X}}}+\tilde{\mathcal{D}}^h\sim_{\mathbb Q} \tilde{\pi}^* K_{\mathbb P ^1}.
\]
Therefore, by \autoref{thm: adj_fiber_space}, $(\tilde{\mathcal{X}},\tilde{\mathcal{D}}^h+\tilde{\mathcal{X}}_{Q_1} + \tilde{\mathcal{X}}_{Q_2})$ is a log Calabi--Yau pair, where $Q_1$ and $Q_2$ are distinct closed points of $\mathbb P ^1$.

Now, we set
\[
    U \coloneqq \big( \tilde{\mathcal{X}} \setminus (\tilde{\mathcal{D}}^h+\tilde{\mathcal{X}}_{Q_1}+\tilde{\mathcal{X}}_{Q_2}) \big)_{\rm reg} .
\]
By construction, $U$ is smooth and satisfies $\overline{\kappa}(U)=0$.
We denote by $g \colon U \to \mathbb G _m$ the morphism induced by $\tilde \pi$.
By \autoref{sec:key_example}, two general fibers of $g$ are not birational.

We claim that $g$ is the quasi-Albanese morphism of $U$.
To show this claim, we first argue that the quasi-Albanese morphism of a general fiber of $g$ is trivial.
A general fiber of $g$ has the form $\tilde{\mathcal{X}}_P \setminus \tilde{\mathcal{D}}^h_P$, where $\tilde{\mathcal{X}}_P $ is a degree $d$ hypersurface of dimension $(d-1)$ and $\tilde{\mathcal{D}}^h_P$ is a smooth anti-canonical $K$-trivial hypersurface in $\tilde{\mathcal{X}}$.
We consider the short exact sequence
\begin{equation*}
    0 \to \Omega_{\tilde{\mathcal{X}}_P}^1 \to \Omega_{\tilde{\mathcal{X}}_P}^1(\log \tilde{\mathcal{D}}^h_P) \to \mathcal{O}_{\tilde{\mathcal{D}}^h_P} \to 0.
\end{equation*}
Since $\tilde{\mathcal{X}}_P$ is a Fano variety, we have $H^0 \big( \tilde{\mathcal{X}}_P,\Omega_{\tilde{\mathcal{X}}_P}^1 \big) = 0$.
In turn, the morphism 
$$
    H^0 \big( \tilde{\mathcal{X}}_P,\Omega_{\tilde{\mathcal{X}}_P}^1(\log \tilde{\mathcal{D}}^h_P) \big) \to H^0 \big( \tilde{D}^h_P , \mathcal{O}_{\tilde{\mathcal{D}}^h_P} \big)
$$ 
is injective.
Then, to deduce that 
$H^0 \big( \tilde{\mathcal{X}}_P,\Omega_{\tilde{\mathcal{X}}_P}^1(\log \tilde{\mathcal{D}}^h_P) \big) = 0$, it suffices to show that the connecting homomorphism 
$\delta \colon H^0 \big( \tilde{\mathcal{D}}^h_P,\mathcal{O}_{\tilde{\mathcal{D}}^h_P} \big) \to H^1 \big( \tilde{\mathcal{X}}_P,\Omega_{\tilde{\mathcal{X}}_P}^1 \big)$ is injective.
But $\delta(1)=c_1(\tilde{\mathcal{D}}^h_P)$ holds by \cite[Theorem 12.2]{GKKP}, and thus $\delta$ is injective.

By the construction in \cite[\S~3]{Fuj24}, 
$h^0 \big( \tilde{\mathcal{X}}_P,\Omega_{\tilde{\mathcal{X}}_P}^1(\log \tilde{\mathcal{D}}^h_P) \big) = 0$ 
computes the dimension of the quasi-Albanese variety of $\tilde{\mathcal{X}}_P \setminus \tilde{\mathcal{D}}^h_P$.
Thus, we conclude that the quasi-Albanese morphism of a general fiber of $g$ is trivial.

Denote by $\alb_U \colon U \to V$ the quasi-Albanese morphism of $U$.
By \cite[Theorem 1.2]{Fuj24}, $\alb_U$ is dominant with irreducible general fiber.
As a general fiber of $g$ has trivial quasi-Albanese, it is contracted by the quasi-Albanese morphism of $U$ by the universal property of the quasi-Albanese morphism.
In turn, for dimension reasons, $V$ has dimension 0 or 1.
Since $g$ is a non-trivial morphism to a torus, we deduce that $V = \mathbb G _m$.
Since both $g$ and $\alb_U$ have irreducible general fiber, we conclude that $\alb_U = g$.
In particular, two general fibers of the quasi-Albanese morphism of $U$ are not birational.
\end{proof}

\subsection{An example without boundary} 
\label{sec:no--boundary}

In this subsection we aim at birationally modifying the example in \autoref{prop: plt-counterexample} to obtain an example without boundary. To this end, we first study the generic fiber $(\mathcal{Y}_\eta, \mathcal{B}_\eta)$ of the example in \autoref{prop: plt-counterexample}.

\begin{lemma} 
    \label{lem: crepant_model_nobound}
    There exists a crepant birational model $(\mathcal{Y}''_{\eta}, 0)$ of $(\mathcal{Y}_\eta, \mathcal{B}_\eta)$.
\end{lemma}
 
\begin{proof}
    For ease of notation, we write $Y \coloneqq \mathcal{Y}_\eta$, which is a smooth degree $d$ hypersurface in $\mathbb P ^d_{\eta}$, and $S \coloneqq \mathcal{B}_{\eta}$, which is a smooth $K$-trivial variety. Note that $(Y,S)$ is a log Calabi--Yau pair. Moreover, we have
    \begin{equation} 
        \label{eq:S_S}
        S|_S \sim H|_S ,
    \end{equation}
    where we regard $S$ as a divisor in $Y$ and $H$ as an hyperplane section in $\mathbb P ^d$.
    We now fix a general element $\Gamma \in \big| 2H|_S \big|$.
    
    Denote by $Y'$ the blow-up of $Y$ along $\Gamma$, and by $\pi$ the corresponding birational morphism.
    We then have $\pi^*S = S'+F'$, where $S'$ denotes the strict transform of $S$ and $F'$ denotes the reduced $\pi$-exceptional divisor.
    
    Since $\Gamma \subset S$, we infer that $\pi|_{S'} \colon S' \to S$ is an isomorphism.
    As a consequence, we have
    \begin{equation} 
        \label{eq:S_S2}
        S|_S \sim (\pi^*S)|_{S'} \sim F'|_{S'}+S'|_{S'} \sim \Gamma' + S'|_{S'},
    \end{equation}
    where $\Gamma'$ denotes the preimage of $\Gamma$ under the isomorphism $\pi|_{S'}$.
    Since by \eqref{eq:S_S} we have $2S|_S \sim \Gamma$, from \eqref{eq:S_S2} we obtain
    \begin{equation} \label{eq:S'_S'}
        S'|_{S'} \sim - S|_{S}.
    \end{equation}
    
    Now, we consider the following short exact sequence
    \begin{equation}
        0 \to \mathcal{O}_{Y'}(\pi^*S) \to \mathcal{O}_{Y'}(\pi^*S+S') \to \mathcal{O}_{S'} \big( (\pi^*S+S')|_{S'} \big) \to 0.
    \end{equation}
    By \eqref{eq:S_S}, $|\mathcal{O}_{Y'}(\pi^*S)|$ is a free linear series.
    Thus, the base locus of $|\mathcal{O}_{Y'}(\pi^*S+S')|$ is contained in $S'$.
    By \eqref{eq:S_S2} and \eqref{eq:S'_S'}, we have 
    \[ 
        (\pi^*S+S')|_{S'} \sim 2S|_S+2S'|_S \sim 0 .
    \]
    Thus, $\mathcal{O}_{S'} \big( (\pi^*S+S')|_{S'} \big)$ is globally generated.
    
    Then, in order to conclude that $|\mathcal{O}_{Y'}(\pi^*S+S')|$ is free, it suffices to show that $\Gamma \big(Y',\mathcal{O}_{Y'}(\pi^*S+S') \big) \to \Gamma \big(S',\mathcal{O}_{S'}((\pi^*S+S')|_{S'}) \big)$ is surjective.
    In turn, this follows from the vanishing of $H^1 \big(Y',\mathcal{O}_{Y'}(\pi^*S) \big)$.
    This vanishing holds by the following chain of equalities
    \[
    H^1 \big(Y',\mathcal{O}_{Y'}(\pi^*S) \big) = H^1 \big(Y,\mathcal{O}_{Y}(S) \big) = 0,
    \]
    where the first equality follows from the projection formula and the fact that $Y$ has rational singularities, and the second equality follows from the Kodaira vanishing theorem, since $-K_Y$ is ample.
    
    Therefore, the divisor $\pi^*S+S'$ is big and semi-ample.
    Furthermore, by construction, it is trivial along $S'$.
    Hence, a sufficiently large multiple of $\pi^*S+S'$ defines a birational morphism $\psi \colon Y' \to Y''$ to a normal variety that contracts the divisor $S'$ to a point.
    
    Now, we consider the log Calabi--Yau pair $(Y,S)$, which is log smooth, 
    and we note that $\pi^*(K_Y+S)=K_{Y'}+S'$, since $\Gamma \subset S$ holds.
    In turn, since $K_Y+S\sim 0$ and since $\psi$ contracts $S'$, we deduce that $Y''$ is a log canonical variety with $K_{Y''} \sim 0$.
    In particular, $(Y'',0)$ is crepant birational to $(Y,S)$.
\end{proof}

We are ready to prove \autoref{thm: counterexample_lc}.
We replace the generic fiber of $(\mathcal{Y}, \mathcal{B})$ with the model constructed in \autoref{lem: crepant_model_nobound}, and after spreading out and performing some birational modification, we can produce an example without boundary.

\begin{proof}[Proof of \autoref{thm: counterexample_lc}]
    Let $\Gamma \subset \mathcal{Y}_{\eta}$ be as in the proof of \autoref{lem: crepant_model_nobound}, and let $\mathcal{C}$ denote its closure in $\mathcal{Y}$.
    We choose a nonempty open subset $V \subset E$ such that $(\mathcal{Y}|_V, \mathcal{B}) \to V$ is log smooth over the base and $\mathcal{C}|_V \to V$ is also smooth.
    
    We construct the claimed model by performing crepant birational modifications on $(\mathcal{Y}, \mathcal{B})$ in two steps.
    
    \medskip
    
    {\bf Step 1}: We perform blow-ups along $\mathcal{B}$ to make it more negative.
    
    First, we blow up $\mathcal{C} \subset \mathcal{B}$ with $\psi \colon \mathrm{Bl}_{\mathcal{C}}\mathcal{Y} \to \mathcal{Y}$.
    Then, we take a thrifty log resolution of the pair $(\mathrm{Bl}_{\mathcal{C}}\mathcal{Y}, \psi_{*}^{-1}\mathcal{B}+F)$, where $F$ are the singular fibers of $\alb_{\mathrm{Bl}_{\mathcal{C}}\mathcal{Y}}$.
    We denote the final model by $\varphi \colon \mathcal{Y}' \to \mathcal{Y}$ and we observe that the generic fiber $\mathcal{Y}'_{\eta}$ corresponds to the intermediate model $Y'$ in the proof of \autoref{lem: crepant_model_nobound}.
    
    \medskip
    
    {\bf Step 2}: We run two MMPs to contract $\mathcal{B'}$, the strict transform of $\mathcal{B}$ on $\mathcal{Y}'$.
    
    Denote by $\mathcal{E}'$ the sum of the $\varphi$-exceptional divisors that are vertical over $E$.
    By construction, for $0 < \epsilon \ll 1$, the pair $(\mathcal{Y}',\mathcal{B'}+\epsilon \mathcal{E}')$ is plt.
    Over $V$, this pair is crepant birational to $(\mathcal{Y},\mathcal{B})$.
    Furthermore, since $(\mathcal{Y},\mathcal{B})$ is plt and $K_{\mathcal{Y}}+\mathcal{B}$ is Cartier, the pair $(\mathcal{Y}, \mathcal{B})$ is actually canonical, i.e., every $\varphi$-exceptional divisor has discrepancy at least 0 with respect to $(\mathcal{Y},\mathcal{B})$.
    Thus, by \cite[Lemma 2.10]{Lai11}, $\supp (\mathcal{E}')$ is in the relative diminished base locus of $(\mathcal{Y}',\mathcal{B'}+\epsilon \mathcal{E}')$ over $E$.
    By \cite[Theorem 1.1]{HX13}, we can run a $(K_{\mathcal{Y}'}+\mathcal{B'}+\epsilon \mathcal{E}')$-MMP over $E$, which ends with a relative good minimal model $(\mathcal{Y}'',\mathcal{B}'')$. 
    By construction, this minimal model is isomorphic to $(\mathcal{Y}',\mathcal{B'}+\epsilon \mathcal{E}')$ over the generic point of $E$ and contracts exactly $\mathcal{E}'$.
    In particular, $\mathcal{Y}'' \dashrightarrow \mathcal{Y}$ is a rational contraction that contracts only one divisor that is horizontal over $E$, which we denote by $\mathcal{F}''$.
    Furthermore, $(\mathcal{Y}'',\mathcal{B}'')$ is crepant birational to $(\mathcal{Y},\mathcal{B})$.
    In particular, $(\mathcal{Y}'',\mathcal{B}'')$ is plt.
    
    By the construction in \autoref{lem: crepant_model_nobound}, the restriction of $2\mathcal{B}''+\mathcal{F}''$ to the generic fiber $\mathcal{Y}''_\eta$ is big and semi-ample.
    Thus, we may find $0 \leq \mathcal{G}'' \sim_{\Q,g''} 2\mathcal{B}''+\mathcal{F}''$ whose support does not contain $\mathcal{B}''$.
    Then, for $0 < \epsilon \ll 1$, the pair $(\mathcal{Y}'',\mathcal{B}''+\epsilon \mathcal{G}'')$ is still plt, and hence by \cite[Theorem 1.1]{HX13}, $(\mathcal{Y}'',\mathcal{B}''+\epsilon \mathcal{G}'')$ admits a relative good minimal model over $E$.
    We denote this model by $\mathcal{Y}'''$.
    By \autoref{lem: crepant_model_nobound}, $\mathcal{Y}'' \dashrightarrow \mathcal{Y}'''$ contracts $\mathcal{B}''$.
    Thus, since $(\mathcal{Y}'',\mathcal{B}'')$ is log Calabi--Yau, we conclude that $\mathcal{Y}'''$ is crepant birational to $(\mathcal{Y}'',\mathcal{B}'')$.
    In particular, $X \coloneqq \mathcal{Y}'''$ is the crepant birational model of $(\mathcal{Y},\mathcal{B})$ that is claimed in the statement.
    By \autoref{prop: plt-counterexample} we have that any two general fibers of $\Alb_X$ are not birational to each other, thus showing \autoref{item: not_bir_fiber}.
    
    The statement \autoref{item: univ_cover} on the universal cover follows as in \autoref{prop: univ_cover_boundary}.
\end{proof}

\section{Counterexamples for Fano varieties}

In this section, we show that the hypotheses on both the positivity of the anti-canonical divisor and the singularities in \autoref{thm: lc_alb_contraction} are necessary.

First, we show that we cannot relax the nefness hypothesis: we give below an example of a smooth projective surface with big but not nef anti-canonical class whose Albanese morphism is not surjective.

\begin{example} 
    \label{example:big_non-nef_nonsurjAlb}
    Let $C$ be a smooth projective curve of genus $g \geq 2$. Pick an integer $m \geq 2$ and set 
    $X \coloneqq \mathbb{P}_C \big( \mathcal{O}_C \oplus \mathcal{O}_C(mK_C) \big)$.
    By the Euler sequence, we have
    \begin{equation} 
        \label{eq:canonical_proj_bunlde}
        K_X = \pi^*(K_C - mK_C) - 2E,
    \end{equation}
    where $E$ is the negative section.
    Furthermore, 
    \begin{equation} 
        \label{eq:pseff_cone}
    \textup{Pseff}(X)= \mathbb {R}_+ [F] + \mathbb{R}_{+}[E],
    \end{equation}
    where $F$ denotes a fiber of the $\mathbb{P}^1$-bundle $\pi \colon X \to C$.
    By \autoref{eq:canonical_proj_bunlde} and \autoref{eq:pseff_cone} we infer that $-K_X$ is big.
    Note also that $X$ is not rationally chain connected.
    Finally, by \autoref{lem:Alb_RCC_fibers} we obtain $\Alb_X \cong \Alb_C$, and since $\alb_C \colon C \to \Alb_C$ is not surjective by the hypothesis $g \geq 2$, we conclude that $\alb_X \colon X \to \Alb_X$ is not surjective either.
\end{example}

Next, we construct a Fano variety of dimension $3$ with worse-than-log canonical singularities and with non-surjective Albanese morphism, proving thus \autoref{thm: counterexample_non_lc}.
This result, together with the previous example, shows that, under the sole assumption that $-(K_X+\Delta)$ is big, \autoref{lem:induced_map_alb} is optimal.
We refer to \cite[\S~3.1]{kk-singbook} for more information about the cone construction that will be used in \autoref{example:nonRCC_verysingular_3fold_nonsurjAlb}.

\begin{example}
    \label{example:nonRCC_verysingular_3fold_nonsurjAlb}
    Let $C$ be a smooth projective curve of genus $g \geq 2$. Consider the smooth projective surface $S  \coloneqq  C \times C$, denote by $\proj_i \colon S \to C$ the projection to the $i$-th factor, 
    and note that $\omega_S = \proj_1^* \omega_C \otimes \proj_2^* \omega_C$ is ample.
    Fix a positive integer $n \geq 1$, and set
    \[
    \cM \coloneqq \omega_S^{\otimes n} \otimes \big( \cO_C \boxtimes \omega_C \big) = \omega_S^{\otimes n} \otimes \proj_2^* \omega_C .
    \]
    Note that $\cM = \cO_S(M)$ is an ample line bundle on $S$, and we have
    \begin{equation}
        \label{eq:ample_div}
        M \sim nK_S + \pr_2^* K_C .
    \end{equation}
    Consider the projective cone $Z \coloneqq  C_p(S,\cM)$ over $S$ with conormal bundle $\cM$, the blow-up $ \rho \colon Y \to Z $ of its vertex, where
    \[ 
    Y \coloneqq  BC_p(S,\cM) = \mathbb{P}_S (\cO_S \oplus \cM) , 
    \]
    and denote by $E \cong S$ the $\rho$-exceptional divisor on $Y$; see \cite[\S~3.1]{kk-singbook}. 
    \begin{center}
        \begin{tikzcd}[row sep = large, column sep = large]
            Y \arrow[r, "\pi"] \arrow[d, "\rho" swap] & S \arrow[d, "\pr_2"] \\
            Z & C
        \end{tikzcd}
    \end{center}
    Note that $Z$ is a normal projective variety such that $K_Z$ is not $\Q$-Cartier, see \cite[Proposition 3.14]{kk-singbook}, and that $Y$ is a smooth projective variety satisfying 
    \begin{equation}\label{eq:can_blow_up_vertex}
        K_Y \sim \pi^* \big( (1-n) K_S - \pr_2^* K_C \big) - 2E .
    \end{equation}
    Furthermore, we have
    \begin{equation}\label{eq:E_E}
        E|_E \sim -M,
    \end{equation}
    cf.\ \cite[Proposition V.2.8]{Har77}.
    
    Let $H_Z$ be the hyperplane section on $Z$ and consider the Cartier divisor
    \[
    L  \coloneqq  (n-1)\rho^* H_Z + \pi^* \pr_2^* K_C
    \]
    on $Y$.
    Since $L$ is big and semi-ample, it determines a birational fibration to a normal variety $g \colon Y \to X$ such that $ L \sim_{\bQ} g^* A_X$, where $A_X$ is an ample $\bQ$-Cartier $\bQ$-divisor on $X$.
    The morphism $g$ contracts exactly the curves $\gamma \subseteq Y$ such that $L \cdot \gamma=0$.
    As both $H_Z$ and $K_C$ are nef, this means that the contracted curves satisfy $\rho^*H_Z \cdot \gamma=0$ and $\pi^* \pr_2^* K_C \cdot \gamma=0$.
    The first condition implies that $\gamma$ is in the exceptional locus of $\rho$, and hence in $E$.
    Out of such curves, only the fibers of $(\pi \circ \pr_2)|_E$ intersect $\pi^* \pr_2^* K_C$ trivially, so they are the only curves on $Y$ contracted by $g$.
    Observe that, under the natural isomorphism $E \cong S$, the morphism $(\pi \circ \pr_2)|_E$ coincides with $\pr_2$.
    
    It follows from the Rigidity Lemma \cite[Proposition 1.14]{Deb01} that there exists a unique morphism $u \colon X \to Z$ such that $ \rho = u \circ g $ and a unique morphism $v \colon X \to C$ such that $\pr_2 \circ \pi = v \circ g$.
    \begin{center}
        \begin{tikzcd}[row sep = large, column sep = large]
            Y \arrow[rr, "\pi"] \arrow[dr, "g"] \arrow[dd, "\rho" swap] && S \arrow[dd, "\pr_2"] \\
            & X \arrow[dl, "u"] \arrow[dr, "v" swap] \\
            Z && C.
        \end{tikzcd}
    \end{center}
    
    Fix a point $c \in C$. Then $\pi^{-1} \left( \pr_2^{-1} (c) \right)$ is a ruled surface over $C$, and $v^{-1}(c)$ is obtained from it by contracting a section.
    In particular, the fiber $v^{-1}(c)$ is rationally chain connected.
    By \autoref{lem:Alb_RCC_fibers}, we obtain $\Alb_X \cong \Alb_C$, and since $\alb_C \colon C \to \Alb_C$ is not surjective, $\alb_X \colon X \to \alb_X$ is not surjective either. Therefore, the normal projective variety $X$ is not rationally chain connected by \autoref{lem:RCC_surj_Alb}.

    We will now show that $X$ has worse-than-log canonical singularities and ample anti-canonical class.
    To this end, we first show that $\rho(Y/X)=1$.
    Since $Y$ is a projective bundle over $S$, its N\'eron--Severi group is generated by the N\'eron--Severi group of $S$ and any $\pi$-ample Cartier divisor.
    Since $\rho^*H_Z$ is trivial along $E$, it follows that $\rho^*H_Z$ is $\pi$-ample.
    Since any two fibers of $E \to C$ are numerically equivalent as curves in $E$ and any curve in $E$ has zero intersection with $\rho^*H_Z$, we deduce that any two fibers of $E \to C$ are numerically equivalent as curves in $Y$.
    In particular, $\rho(Y/X) = 1$ holds.
    In turn, we have $K_Y + aE \equiv_g 0$ for some $a \in \Q$.
    
    To show that $K_X$ is $\Q$-Cartier, it suffices to show that $K_Y + aE \sim_{\Q,g} 0$.
    Fix a curve $\gamma$ that is contracted by $\rho$ and notice that $\pi_*\gamma \neq 0$.
    In view of the above construction, by \eqref{eq:ample_div}, \eqref{eq:can_blow_up_vertex}, and \eqref{eq:E_E} we obtain
    \begin{align*}
        0 &= (K_Y + aE) \cdot \gamma \\
        &= (a-2) (E \cdot \gamma) + (1-n) (\pi^* K_S \cdot \gamma) \\ 
        &= -(a-2) (M \cdot \pi_* \gamma) + (1-n) (K_S \cdot \pi_* \gamma) \\
        &= -n(a-2) (K_S \cdot \pi_* \gamma) + (1-n) (K_S \cdot \pi_* \gamma) ,
    \end{align*}
    which yields 
    \[
    a = 1 + \frac{1}{n} ,
    \]
    since we have $K_S \cdot \pi_*\gamma>0$ by the ampleness of $K_S$.
    Set $D \coloneqq K_Y + \big( 1 + \frac{1}{n} \big) E$ and $H_Y \coloneqq \rho^* H_Z$. By adjunction on $E$ or $H_Y$ we infer that $K_Y = - E - H_Y + \pi^* K_S$, and thus $H_Y = E + \pi^* M$ by \eqref{eq:can_blow_up_vertex}. Therefore, we have
    \begin{align*}
        -D &= -K_Y - \left(1 + \frac{1}{n} \right) E = E+H_Y-\pi^*K_S-\left(1+\frac{1}{n}\right)E \\[0.25em]
        &= \frac{1}{n} \pi^*M-\frac{1}{n}  H_Y+H_Y-\pi^*K_S = \frac{1}{n} \big( \pi^* \operatorname{pr}_2^* K_C + (n-1)H_Y \big) \\[0.25em]
        &= \frac{1}{n}L \sim_\Q \frac{1}{n} g^* A_X ,
    \end{align*}
    which yields that $-K_X \sim_\Q \frac{1}{n} A_X$ is an ample $\Q$-Cartier divisor.
    
    Finally, the relation
    \[
    K_Y \sim_\Q g^* K_X - \left( 1 + \frac{1}{n} \right) E
    \]
    shows that the singularities of $X$ are worse-than-log canonical, as asserted.
\end{example}

\appendix

\section{On the polystability of the tangent sheaf -- by Niklas M\"uller}
  
    One of the deepest insights into the geometry of smooth $K$-trivial varieties is the existence of K\"ahler--Einstein metrics established by Yau \cite{Yau_CalabiConjecture}. Two key consequences are the polystability of their tangent bundle \cite{Luebke_StabilityOfHEBundles, Kobayashi_StabilityOfHEBundles} and the Beauville--Bogomolov decomposition \cite{Beau83}. During the last years, the existence of K\"ahler--Einstein metrics \cite{EGZ_SingularKEMetrics}, the polystability of the tangent bundle \cite{guenancia_SemistabilityOfTangent} and the decomposition theorem \cite{HP19, BGL_DecompositionOfKahlerCYs} have all been extended to $K$-trivial varieties with klt singularities. 
    
    In the main text, Bernasconi--Filipazzi--Patakfalvi--Tsakanikas construct $K$-trivial varieties $X$ with log canonical singularities for which the Albanese morphism $\alpha\colon X \rightarrow \Alb_X$ is not (birationally) isotrivial. This implies that the analogue of the Beauville--Bogomolov decomposition theorem does not hold for such varieties. In this appendix we strengthen their results by showing that the tangent sheaf of $X$ fails to be polystable. Moreover, we show that $\alpha\colon X \rightarrow \Alb_X$ need not even be locally topological trivial, answering a question of Xu \cite{Xu_kltCYs} in the process.
    
    \begin{theorem}\emph{(c.f.\ \autoref{thm: counterexample_lc})}
        For any integer $d\geq 4$, there exists a projective log canonical variety $X$ of dimension $d$ such that\label{thm:IntroExampleNoBoundary}
        \begin{itemize}
            \item[(1)] the canonical divisor $K_X$ is Cartier and $K_X\sim 0$;
            \item[(2)] the tangent sheaf $\mathcal{T}_{X}$ is not polystable with respect to any polarisation;
            \item[(3)] the Albanese morphism $\alpha\colon X \rightarrow \Alb_X$ is not a topological fibre bundle.
        \end{itemize}
    \end{theorem}
    In contrast, the tangent sheaf of a variety $X$ with ample canonical divisor is polystable even when $X$ has log canonical singularities \cite{guenancia_SemistabilityOfTangent}. We point out that our proof of \autoref{thm:IntroExampleNoBoundary}.(3) shows that $\alpha$ need not even be a $\Q$-homology fibre bundle. 

    Finally, adopting slightly the construction in \autoref{ex: 3-dimensional}, we can also give a log smooth counter-example:
    \begin{proposition}
        There exists a projective log smooth threefold pair $(X, D)$ such that\label{Prop:IntroLogSmoothEx}
        \begin{itemize}
            \item[(1)] the log canonical divisor satisfies $K_{X} + D \sim 0$;
            \item[(2)] the logarithmic tangent bundle $\mathcal{T}_{X}(-\log D)$ is not polystable with respect to any polarisation; 
            \item[(3)] the (quasi-) Albanese morphism $\alpha\colon X\setminus D \rightarrow \Alb_{X\setminus D}$ is not a topological fibre bundle.
        \end{itemize}
    \end{proposition}
    In particular, the Albanese morphism $\alpha\colon (X, D) \rightarrow \Alb_X$ is not a topological fibre bundle of pairs. 
    
    Note that, since the appearance of this article, Collins--Guenancia \cite{CG25} have constructed further examples satisfying the conclusion in \autoref{Prop:IntroLogSmoothEx} but with the additional property that the analogue of the Bochner property fails on $X$ and $X\setminus D$ admits a complete Ricci flat K\"ahler metric.
    
    \subsubsection*{\normalfont \textbf{Proof of Proposition A.2}}
    \label{subsec:LogSmoothExample}
    
    The following construction is heavily inspired by the proof of \autoref{ex: 3-dimensional}. Our construction has the advantage that the example in \autoref{ex: 3-dimensional} does not satisfy property $(3)$ in \autoref{Prop:IntroLogSmoothEx}.
    \begin{construction} \label{const}
        Let $p_1, p_2, q_1, q_2 \in \mathbb{P}^2$ be four general points. Let $\ell_i$ be the unique line through $p_i, q_i$ for $i= 1,2$. Choose a general line $\ell_3$, meeting the line $\ell_i$ in the single point $r_i$, respectively.
        Let $\pi\colon X \rightarrow \mathbb{P}^2$ be the blow-up of the $4$ points $p_1, p_2, q_1, q_2$. Let $L_i$ denote the strict transform of $\ell_i$, $i = 1,2,3$. Then, setting  $D := L_1 + L_2 + L_3$, one readily checks that
        $$
        K_X + D 
        \sim \pi^*(K_{\mathbb{P}^2} +\ell_1 + \ell_2 + \ell_3)
        \sim 0.
        $$
        Now, let $\phi\colon E \rightarrow L_3 \cong \mathbb{P}^1$ be a finite branched cover, unramified over $r_1, r_2$, and let us denote $\phi^{-1}(r_i) = \{t_i, s_i\}$, $i=1,2$. Set
        $$
        \mathcal{X} := X \times E,\quad
        \mathcal{D} := D \times E,\quad
        \mathcal{L}_i := L_i\times E
        $$
        and let
        $$
        \Gamma = \big\{ (\phi(t), t)\ \big| \ t\in E \big\} \subseteq L_3 \times E \subseteq \mathcal{D} \subseteq  \mathcal{X}
        $$
        be the graph of $\phi$. Then $\Gamma \cong E$ via the second projection.

        We consider the blow-up $\psi\colon \mathcal{Y} \rightarrow \mathcal{X}$ along $\Gamma$. Let $\mathcal{B} := \widehat{\mathcal{L}}_1 + \widehat{\mathcal{L}}_2 + \widehat{\mathcal{L}}_3$ be the strict transform of $\mathcal{D}$ under $\psi$ and let $\mathcal{E}$ be the exceptional divisor of $\psi$. Then
        $$
        K_{\mathcal{Y}} \sim \psi^*K_{\mathcal{X}} + \mathcal{E},\quad
        \psi^*\mathcal{L}_1 \sim \widehat{\mathcal{L}}_1,\quad
        \psi^*\mathcal{L}_2 \sim \widehat{\mathcal{L}}_2,\quad
        \psi^*\mathcal{L}_3 \sim \widehat{\mathcal{L}}_3 + \mathcal{E}.
        $$
        In particular,
        $$
        K_{\mathcal{Y}} + \mathcal{B} \sim
        \psi^*\Big(K_{\mathcal{X}} + \mathcal{D}\Big) \sim 0.
        $$
        Note that $\psi$ restricts to an isomorphism $\widehat{\mathcal{L}}_3 \cong \mathcal{L}_3$, while for $i=1,2$ the natural map $\widehat{\mathcal{L}}_i \rightarrow \mathcal{L}_i$ coincides with the blow-up of $\mathcal{L}_i$ in the two points $(r_i, t_i), (r_i, s_i)$.

        Now, consider the morphism
        \begin{align}
          \alpha:= \psi \circ \textmd{pr}_2\colon \mathcal{Y} \rightarrow E,  \label{eqs-ExLogSmoothExample-Albanese}
        \end{align}
        which is clearly the Albanese morphism of $\mathcal{Y}$.
        Then $Y_t := \alpha^{-1}(t)$ is simply the blow-up $\varphi_t\colon Y_t\rightarrow \mathbb{P}^2$ of the $5$ distinct points $p_1, p_2, q_1, q_2, \phi(t)$. Moreover,
        $B_t := \mathcal{B}\cap Y_t$ can be described as follows: For 
        $t\in E\setminus \{t_1, t_2, s_1, s_2\}$, the divisor $B_t = \widehat{L}_1 + \widehat{L}_2 + \widehat{L}_3$ is the strict transform of $D$ and in fact isomorphic to $\ell_1+\ell_2+\ell_3$. In particular, $B_t$ is a closed chain of three rational curves. On the other hand, for the special values $t \in \{t_1, t_2, s_1, s_2\}$ we have that $B_t = \widehat{L}_1 + \widehat{L}_2 + \widehat{L}_3 + \mathcal{E}_t$ is a closed chain of four rational curves.

        \begin{figure*}[h!]
            \centering
            \includegraphics[width=12.5cm]{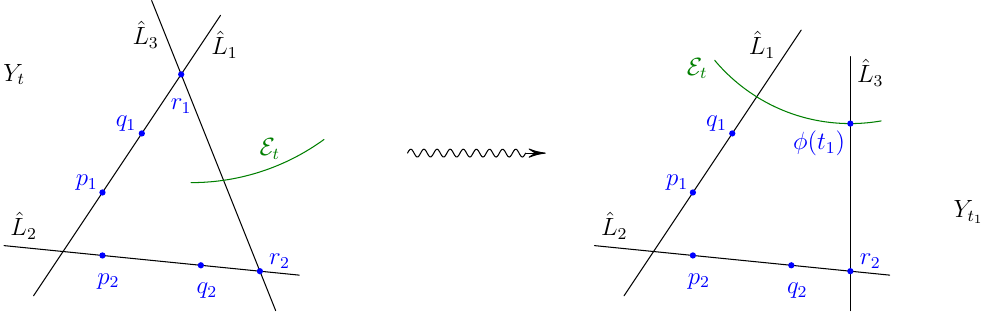}
            \label{fig:enter-label}
            \caption{A general fibre $Y_t$ degenerates to the special fibre $Y_{t_1}$.}
        \end{figure*}
    \label{constr:LogSmoothExample}
    \end{construction}
    
    \begin{proposition}
        In the set-up of Construction \ref{constr:LogSmoothExample},
        \begin{itemize}
            \item[(1)] no two general fibres of $\alpha$ are isomorphic. 
            \item[(2)] $\Aut^0(\mathcal{Y}) = \{\emph{id\}}$. 
        \end{itemize}
        \label{PropLogSmoothExampleBasicProperties}
    \end{proposition}
    \begin{proof}
        Assertion $(1)$ can be proved exactly as \autoref{ex: 2dim}.(c). Note that the argument given there also shows that $\Aut^0(Y_t) = \textmd{\{id\}}$ for general $t\in E$.

        Concerning $(2)$, since the Albanese morphism is universal, $\Aut^0(\mathcal{Y})$ acts equivariantly on $\alpha$. In other words, $\Aut^0(\mathcal{Y})$ maps fibres of $\alpha$ onto fibres of $\alpha$. From $(1)$, it follows that this action has to map each (general) fibre of $\alpha$ onto itself. In particular, $\Aut^0(\mathcal{Y})$ acts effectively on a general fibre $Y_t$ of $\alpha$. However, from the proof of $(1)$ we know that $\Aut^0(Y_t) = \{\textmd{id\}}$. The proposition is proved.
    \end{proof}
    \begin{remark}
        The morphism $\alpha' \colon\mathcal{Y}\setminus\mathcal{B} \rightarrow E$ is the quasi-Albanese of the open variety $\mathcal{Y}\setminus\mathcal{B}$. Indeed, as in the proof of \autoref{thm: non-isotriviality-quasialbanese}, it suffices to show that the quasi-Albanese of any fibre of $\alpha'$ is trivial. By \cite[Proof of Lemma 2.3]{Fuj24} we have
        $$
        H^0\Big(Y_t, \Omega_{Y_t}^1(\log B_t)\Big)
        \subseteq  H^0\Big(Y_t, \Omega_{Y_t}^1\Big(\log \varphi_t^{-1}\left(\bar{D}\right)\Big)\Big)
        = H^0\Big(\mathbb{P}^2, \Omega_{\mathbb{P}^2}^1\left(\log \bar{D}\right)\Big),
        $$
        where $\varphi_t\colon Y_t\rightarrow \mathbb{P}^2$ and $\bar{D} = \ell_1  + \ell_2 + \ell_3 \subseteq \mathbb{P}^2$. Now $\Omega_{\mathbb{P}^2}^1(\log \bar{D}) \cong \mathcal{O}_{\mathbb{P}^2}^{\oplus 2}$ has two linearly independent global sections. Expressing these forms in local coordinates, one sees that they all have poles along one of the components of $\varphi_t^{-1}(\bar{D})\setminus B_t$.
        \label{rem:QuasiAlbaneseLogSmoothExample}
    \end{remark}
    
    We are now in position to prove that the pair $(\mathcal{Y}, \mathcal{B})$ satisfies the conclusions in Proposition \ref{Prop:IntroLogSmoothEx}:
    \begin{proof}[Proof of Proposition \ref{Prop:IntroLogSmoothEx}.(2)]
        If $\mathcal{T}_{\mathcal{Y}}(-\log \mathcal{B})$ were polystable, then the generically surjective morphism
        $
        \mathcal{T}_{\mathcal{Y}}(-\log \mathcal{B}) \rightarrow  \alpha^*\mathcal{T}_E
        $
        would split. In other words, $\alpha^*\mathcal{T}_E$ would be a direct summand of $\mathcal{T}_{\mathcal{Y}}(-\log \mathcal{B})$. In particular,
        \begin{equation}
        h^0\Big(\mathcal{Y}, \mathcal{T}_{\mathcal{Y}}(-\log \mathcal{B})\Big) 
        \geq h^0\big(\mathcal{Y}, \alpha^*\mathcal{T}_E\big)
        = h^0\big(\mathcal{Y}, \mathcal{O}_{\mathcal{Y}}\big) 
        = 1. \label{Eq:SurjectionAutomorphismGroups}
        \end{equation}
        But as $h^0(\mathcal{Y}, \mathcal{T}_{\mathcal{Y}}(-\log \mathcal{B})) = \dim(\Aut(\mathcal{Y}, \mathcal{B}))$ by
        \cite[Lemma 2.1]{Xu_kltCYs} equation (\ref{Eq:SurjectionAutomorphismGroups}) contradicts Proposition \ref{PropLogSmoothExampleBasicProperties}.(2).
    \end{proof}
    
    \begin{proof}[Proof of \autoref{Prop:IntroLogSmoothEx}.(3)]
        Recall from the \autoref{const} that the numbers $h^*(Y_t, \Q)$ are independent of $t\in E$. Recall, moreover, that $B_t$ is a closed chain of either three or four smooth rational curves, depending on whether or not $t\in\{ t_1, t_2, s_1, s_2\}$. An easy calculation, c.f.\ \cite[Section II.2]{BHPvV_ComactComplexSurfaces}, shows that the singular cohomology of a closed chain $C_r$ of $r$ rational curves satisfies
        $$
        h^*(C_r, \Q) = \left\{
        \begin{array}{cc}
             1, & * = 0 \\
             1, & * = 1 \\
             r, & * = 2 \\
        \end{array}
        \right.
        $$
        In particular, from the long exact sequence
        $$
        \ldots 
        \rightarrow H^*_c(Y_t\setminus B_t, \Q) 
        \rightarrow H^*(Y_t, \Q)
        \rightarrow H^*(B_t, \Q)
        \rightarrow H^{*+1}_c(Y_t\setminus B_t, \Q) 
        \rightarrow \ldots
        $$
        we see that the numbers $h^i_c(Y_t\setminus B_t, \Q)$ are not all independent of $t\in E$. We conclude that $\alpha'\colon \mathcal{Y}\setminus \mathcal{B} \rightarrow E$ is not a topological fibre bundle.
    \end{proof}

    \subsubsection*{\normalfont \textbf{Proof of Theorem A.1}}
	\label{sec:No Boundary}
    
    Let $d\geq 4$ be an integer. In what follows, we show that the log canonical projective $K$-trivial variety of dimension $d$ constructed by Bernasconi--Filipazzi--Patakfalvi--Tsakanikas in the proof of \autoref{thm: counterexample_lc} satisfies the conclusions of \autoref{thm:IntroExampleNoBoundary}. 
    
    We begin by recalling the most important details from their construction: Choose general homogeneous polynomials $f \in \C[x_0, \ldots, x_d]$ of degree $d$ and general homogeneous polynomials $q_1, q_2 \in \C[x_0, \ldots, x_{d-1}]$
    of degree $d-1$. Consider the pencil of hypersurfaces
    $$ 
    \mathcal{X} := V\Big(y_0(f+q_1x_d) + y_1 (f+q_2x_d)\Big) \subseteq \mathbb{P}^d_{[x_0 : \ldots : x_d]} \times \mathbb{P}^1_{[y_0 : y_1]}
    $$
    Set $\mathcal{D} := \mathcal{X} \cap V(x_d)$ and denote by $\pi\colon \mathcal{X}\rightarrow \mathbb{P}^1$ the projection onto the second factor. Let $g\colon E \rightarrow \mathbb{P}^1$ be a finite cover by an elliptic curve, unramified over the finitely many points $y\in \mathbb{P}^1$ for which $(\mathcal{X}_y, \mathcal{D}_y)$ is not log smooth. Denote by
    $\bar{g}\colon (\bar{\mathcal{X}}, \bar{\mathcal{D}}) \rightarrow (\mathcal{X}, \mathcal{D})$
    the normalization of the base-change by $g$. Then $\bar{\pi} = \alpha_{\mathcal{X}} \colon (\bar{\mathcal{X}}, \bar{\mathcal{D}}) \rightarrow E$
    is precisely the Albanese morphism of $\bar{\mathcal{X}}$.
    
    After performing a series of birational modifications, Bernasconi--Filipazzi--Patakfalvi--Tsakanikas arrive at birational models
    $$(\bar{\mathcal{X}}, \bar{\mathcal{D}}) \overset{\Phi}{\dashleftarrow}
    (\mathcal{Y}', \mathcal{B}') \overset{\Psi}{\dashrightarrow}
    \mathcal{Z} := \mathcal{Y}'''$$
    and each of these three pairs is log canonical with $K_{\bar{\mathcal{X}}} + \bar{\mathcal{D}} \sim K_{\mathcal{Y}'} + \mathcal{B}' \sim K_{\mathcal{Z}} \sim 0.
    $
    We are interested in the variety $\mathcal{Z}$, which is precisely the counter-example in \autoref{thm: counterexample_lc}. To study $\mathcal{Z}$, we consider the respective Albanese morphisms:
\[\begin{tikzcd}
	& {\mathcal{Y}'} \\
	{\bar{\mathcal{X}}} && {\mathcal{Z}} \\
	& E
	\arrow["\Phi"', dashed, from=1-2, to=2-1]
	\arrow["\Psi", dashed, from=1-2, to=2-3]
	\arrow["{\alpha_{\mathcal{Y}'}}", from=1-2, to=3-2]
	\arrow["{\alpha_{\bar{\mathcal{X}}}}"', from=2-1, to=3-2]
	\arrow["{\alpha_{\mathcal{Z}}}", from=2-3, to=3-2]
\end{tikzcd}\]
    By construction, there exists a non-empty, Zariski open subset $U\subseteq E$ such that $\alpha_{\bar{\mathcal{X}}}, \alpha_{\mathcal{Y}'}$ and $\alpha_{\mathcal{Z}}$ are topological fiber bundles over $U$ and $\Phi|_U, \Psi|_U$ are everywhere defined and have the following explicit description: Pick $t\in U$ and set 
    $$(X, D) = (\bar{\mathcal{X}}|_t, \bar{\mathcal{D}}|_t) = (\mathcal{X}|_{g(t)}, \mathcal{D}|_{g(t)}),\quad (Y, B) := (\mathcal{Y}', \mathcal{B}'), \quad Z := \mathcal{Z}|_t.$$
    Then $X\subseteq \mathbb{P}^d$ is a Fano hypersurface of degree $d$, $D\subseteq X$ is an anti-canonical hypersurface and $\phi = \Phi_t$ is the blow-up of a smooth divisor $\Gamma = \Gamma_t \in |\mathcal{O}_D(2)|$, depending holomorphically on $t\in U$. Moreover, $B = \hat{D}$ identifies with the strict transform of $D$ under $\phi$ and it is isomorphic to $D$. Finally, $\psi = \Psi_t\colon Y \rightarrow Z$ is contracts $B$ to a point $z_0$ and is an isomorphism everywhere else.

    On the other hand, it seems hard to give a direct description of the fibers of $\alpha_{\mathcal{Y}'}, \alpha_{\mathcal{Z}}$ over $E\setminus U$. Instead, we choose to study these fibers indirectly via the local systems $R^i(\alpha_\bullet)_*\Q|_U$. In what follows, we denote $H^*(-) := H_{\textmd{sing}}^*(-, \Q)$.
    \begin{proposition}
        In the notation above, the natural map of local systems
        \label{identificationsLocalSystems}
        \begin{align}
            \Phi^*\colon R^{d-1}(\alpha_{\bar{\mathcal{X}}})_*\Q|_U
        \hookrightarrow R^{d-1}(\alpha_{\mathcal{Y}'})_*\Q|_U
        \label{eq-local-systems-1}
        \end{align}
        is injective. Moreover, the natural sequence of maps of local systems
        \begin{align}
        R^{d-1}(\alpha_\mathcal{Z})_*\Q|_U
        \overset{\Psi^*}{\longrightarrow} 
        R^{d-1}(\alpha_{\mathcal{Y}'})_*\Q|_U
        \overset{\iota^*}{\longrightarrow}
        R^{d-1}(\alpha_\mathcal{B}')_*\Q|_U 
        \label{eq-local-systems-2}
        \end{align}
        is exact. Furthermore,
        \begin{align}
            R^{d-1}(\alpha_\mathcal{B}')_*\Q|_U  
            \cong R^{d+3}(\alpha_{\mathbb{P}^d\times E})_*\Q
            = H^{d+3}\left(\mathbb{P}^d\right)
            \label{eq-local-systems-3}
        \end{align}
        is isomorphic to the constant local system with fiber $H^{d+3}\left(\mathbb{P}^d\right)$ on $U$.
    \end{proposition}
    \begin{remark}
        By Deligne's semisimplicity theorem \cite[Th\'eor\`eme 4.2.6]{Deligne_TheorieDeHodgeII},  the local system $R^{d-1}(\alpha_{\mathcal{Y}'})_*\Q|_U$ is semisimple. In particular, by (\ref{eq-local-systems-3}), \label{rem:CohomologyOfD}
        $$
        R^{d-1}(\alpha_{\mathcal{Y}'})_*\Q|_U
        \cong \textmd{Im}\big(\Psi^*\big) \oplus \Q,
        \quad \textmd{or}\quad
        R^{d-1}(\alpha_{\mathcal{Y}'})_*\Q|_U
        = \textmd{Im}\big(\Psi^*\big),
        $$
        depending on the parity of $d$ and whether $\iota^*$ is trivial or not.
    \end{remark}
    \begin{proof}[Proof of \autoref{identificationsLocalSystems}]
        Since $\alpha_{\bullet}$ is a topological fiber bundle over $U$, all sheaves in question are local systems. In particular, to verify (\ref{eq-local-systems-1}), it suffices to show that the natural $\pi_1(U)$-equivariant map $\phi^*\colon H^{d-1}(X) \rightarrow H^{d-1}(Y)$ is injective. Since $\phi$ is the blowing-up of a smooth variety along a smooth center, this follows from the well-known formula \cite[Theorem 7.31]{Voisin_HodgeTheoryofAlgebraicVarieties}.

        Regarding (\ref{eq-local-systems-3}), recall that $B\cong D$ via $\phi$. Now, $D\subseteq X \subseteq \mathbb{P}^d$ is a smooth complete intersection of dimension $d-2$. In particular, by Poincar\'e duality and the Lefschetz hyperplane theorem,
        $$
        H^{d-1}(D) \cong H^{d-3}(D)^\vee = H^{d-3}(X)^\vee \cong H^{d-3}\left(\mathbb{P}^d\right)^\vee
        \cong H^{d+3}\left(\mathbb{P}^d\right).
        $$
        Since all of the above identifications are natural in families, they glue to the isomorphism (\ref{eq-local-systems-3}).

        Finally, concerning (\ref{eq-local-systems-2}), we need to verify that the natural sequence of $\pi_1(U)$-equivariant maps
        \begin{align}
            H^{d-1}(Z) \overset{\psi^*}{\rightarrow} H^{d-1}(Y) \overset{\iota^*}{\rightarrow} H^{d-1}(B) \label{eq-local-systems-4}
        \end{align}
        is exact.
        Since $\psi$ is an isomorphism over $Z\setminus \{z_0\}$ and since $\psi^{-1}(z_0) = B$, we see that the $E_2$-page of the Leray spectral sequence $E^{i, j}_2 = H^j(Z, R^i\psi_*\Q) \Rightarrow H^{i+j}(Y, \Q)$ is given by:
        $$
        \begin{array}{cccccccc}
            H^4(B) & 0 & \ldots  \\
            H^3(B) & 0 & 0 & \ldots \\
            H^2(B) & 0 & 0 & 0 & 0 & \ldots \\
            H^1(B) & 0 & 0 & 0 & 0 & 0 & 0 \\
            H^0(Z) & H^1(Z) & H^2(Z) & H^3(Z) & H^4(Z) & H^5(Z) & H^6(Z) & \ldots
        \end{array}
        $$
        In particular, we see that the following sequence is exact:
        $$
        0
        \rightarrow H^{d-2}(B)
        \overset{\delta^{d-2, 0}_{d-1}}{\longrightarrow} H^{d-1}(Z)
        \rightarrow H^{d-1}(Y) 
        \rightarrow H^{d-1}(B)
        \overset{\delta^{d-1, 0}_{d}}{\longrightarrow} H^d(Z)
        \rightarrow 0
        $$
        Here, we denote by $\delta^{i, j}_r \colon E^{i, j}_r \rightarrow E^{i+1-r, j+r}_r$ the differentials. This shows (\ref{eq-local-systems-4}).
    \end{proof}

    \begin{proposition}
        In the notation above it holds that $\Aut^0(\mathcal{Z}) = \{\emph{id}\}$.\label{PropAutExampleNoBoundary}
    \end{proposition}
    \begin{proof}
        Since two general fibers of $\alpha_{\mathcal{Z}}$ are not birational by \autoref{thm: counterexample_lc}, as in the proof of Proposition \ref{PropLogSmoothExampleBasicProperties}.(2) it follows that $\Aut^0(\mathcal{Z})$ maps the fibers of $\alpha_{\mathcal{Z}}$ onto themselves. Hence, $\Aut^0(\mathcal{Z})$ acts effectively on the general fiber $Z$ of $\alpha_{\mathcal{Z}}$. But any automorphism $\mu\in \Aut(Z)$ induces a birational automorphism $\nu \in \textmd{Bir}(X)$. By birational superrigidity of hypersurfaces \cite{Kol19b}, the rational map $\nu \in \Aut(X)$ is biregular. However, $\Aut(X)$ is a finite group by a theorem of Matsumura--Monsky \cite{MatsumuraMonsky_AutomorphismsHypersurfaces}. We conclude that
        $
        \Aut^0(\mathcal{Z}) \subseteq \Aut^0(X) = \{\textmd{id}\}
        $.
    \end{proof}
    
    Xu \cite{Xu_kltCYs} asked whether the quasi-Albanese of a log canonical log $K$-trivial pair is always a homogeneous fibration. \autoref{PropAutExampleNoBoundary} answers this question in the negative.

    Finally, let us prove that the variety $\mathcal{Z}$ satisfies the conclusions in \autoref{thm:IntroExampleNoBoundary}:
    \begin{proof}[Proof of \autoref{thm:IntroExampleNoBoundary}.(2)]
        Ad verbatim as in the proof of Proposition \ref{Prop:IntroLogSmoothEx}.(2).
    \end{proof}
    
    \begin{proof}[Proof of \autoref{thm:IntroExampleNoBoundary}.(3)]
        Heading for a contradiction, we assume that $\alpha_{\mathcal{Z}}$ is a topological fiber bundle. Then $R^{d-1}(\alpha_\mathcal{Z})_*\Q$ is a local system on $E$. In particular, for any $t\in E\setminus U$ and any small analytic disc $t\in \Delta\subseteq E$, the local system $R^{d-1}(\alpha_\mathcal{Z})_*\Q|_{\Delta^\times}$ has trivial monodromy around $t$. By Proposition \ref{identificationsLocalSystems} and Remark \ref{rem:CohomologyOfD}, the same is true of
        $R^{d-1}(\alpha_{\bar{\mathcal{X}}})_*\Q|_{\Delta^\times}$. Here we used that subquotients of the trivial representation are trivial. Recall, that $\bar{\mathcal{X}} \rightarrow E$ was defined as the base change of $\pi\colon \mathcal{X} \rightarrow \mathbb{P}^1$ along a general morphism $g\colon E \rightarrow \mathbb{P}^1$. Since $g$ is \'etale at all points outside $U$, it follows that also the monodromy of
        $$
        R^{d-1}\pi_*\Q|_{\Delta^\times} \cong R^{d-1}(\alpha_{\bar{\mathcal{X}}})_*\Q|_{\Delta^\times}
        $$
        vanishes. Set $V := \mathbb{P}^1 \setminus (g(E\setminus U))$ and let $j\colon V \hookrightarrow \mathbb{P}^1$ be the inclusion. Then the discussion above implies that
        $
        j_*( R^{d-1}\pi_*\Q|_V )
        $
         is a local system on $\mathbb{P}^1$ with fiber $H^{d-1}(X)$. Since $\mathbb{P}^1$ is simply connected, this local system is constant. We deduce that also the local system
        $
        R^{d-1}\pi_*\Q|_V
        $
        is constant. By the Griffiths--Deligne--Schmid theorem on the fixed part \cite[Corollary 13.1.11]{CMSP_PeriodMappingsAndPeriodDomains} we conclude that even the variation of Hodge structures
        $$
        \Big(R^{d-1}\pi_*\Q|_V, \big(R^p\pi_*\Omega^q_{\mathcal{X}/C}|_V\big)_{p+q = d-1}\Big)
        $$
        is constant. But this is absurd as it contradicts Donagi's generic Torelli theorem \cite{Donagi_GenericTorelliProjectiveHypersurfaces} and the fact that $\pi$ is not isotrivial.
    \end{proof}

    \subsection*{Acknowledgments}
    
    I am grateful towards Benjamin Church for making me aware of \cite{Xu_kltCYs} and some early discussions on the topic. Many thanks to Henri Guenancia and Yuji Odaka for answering several of my questions regarding prior known results and to Vasily Rogov for discussions on degenerations of Hodge structures. Last but not least, many thanks to Daniel Greb, Andreas H\"oring and Junyan Cao for their supervision and their constant encouragement. While writing this article, I was supported by the DFG Research Training Group 2553, ``Symmetries and classifying spaces: analytic, arithmetic and derived."

\bibliographystyle{amsalpha}
\bibliography{refs}

\end{document}